\documentclass[12pt]{article}%

\usepackage{graphicx}
\usepackage{amsmath}
\usepackage{mathptmx}      
\usepackage{enumerate}
\usepackage{authblk}
\usepackage[misc]{ifsym}
\usepackage{doi}

\usepackage{hyperref}
\usepackage{textcomp}
\usepackage[T1]{fontenc}
\usepackage{amsmath,amssymb,mathrsfs}
\usepackage{multirow}
\usepackage{subfigure}
\usepackage{amsfonts}
\usepackage{booktabs}
\usepackage{bm}

\DeclareMathAlphabet      {\mathbfit}{OML}{cmm}{b}{it}

\begin{document}

\title{A new mathematical model for dispersion of Rayleigh wave and a machine learning based inversion solver
\thanks{The work is partially supported by the Shenzhen Stable Support Fund for College Researches [No. 20200827173701001 \& No. 20200829143245001] and the Guangdong Fundamental and Applied Research Fund [No. 2019A1515110971].
}
}

\author[1]{Jianxun Yang}
\author[1]{Chen Xu \thanks{Corresponding author:  xuchen@smbu.edu.cn}}
\author[2]{Ye Zhang}
\affil[1]{Shenzhen MSU-BIT University, 518172 Shenzhen, China}
\affil[2]{School of Mathematics and Statistics, Beijing Institute of Technology, 100081 Beijing, China}
\date{}  

\maketitle

\begin{abstract}
In this work, by introducing the seismic impedance tensor we propose a new Rayleigh wave dispersion function in a homogeneous and layered medium of the Earth, which provides an efficient way to compute the dispersion curve -- a relation between the frequencies and the phase velocities. With this newly established forward model, based on the Mixture Density Networks (MDN) we develop a machine learning based inversion approach, named as FW-MDN, for the problem of estimating the S-wave velocity from the dispersion curves. The method FW-MDN deals with the non-uniqueness issue encountered in studies that invert dispersion curves for crust and upper mantle models and attains a satisfactory performance on the dataset with various noise structure. Numerical simulations are performed to show that the FW-MDN possesses the characteristics of easy calculation, efficient computation, and high precision for the model characterization\footnote{Our codes are available at \url{https://github.com/chen-research/A-new-math-model-for-Rayleigh-wave-and-a-MDN-inversion-solver}}.
\end{abstract}

\section{Introduction}
\label{sec:1}
The study of surface waves has recently gained a wide attention in seismic sounding since it plays an essential role in revealing Earth's structure from the shallow near-surface to several hundred kilometers deep into the mantle, which depends on the frequencies and data acquisition configurations. As a kind of surface waves, the speed of the Rayleigh wave is slower than that of the P wave and the S wave, so it is more sensitive with respect to the change of the underground medium. Moreover, the energy attenuation of Rayleigh wave is slower than the geometrical divergence of body waves, which is strong and easy to identify in seismic records. Furthermore, the Rayleigh wave can be used to obtain the information about the S-wave velocity structure in various mediums (e.g., the crust and the lithosphere). Hence, the Rayleigh wave is defined to be the principal wave which bounds between a transverse and a longitudinal wave and propagates along the Earth's surface, and it has been widely applied in the modeling of small-depth geophysical sounding \cite{berg2020shear}, deep seismological studies \cite{li2020surface}, and other scientific areas \cite{mora2018inversion}.

The problem of deriving Rayleigh wave dispersion curves from the crust and upper mantle models is an important item in the study of dispersion characteristics and surface wave sounding technology. The main content in this topic is to establish an implicit equation on the relationship between the Rayleigh wave velocity and the wave frequency and other elastic parameters in the known medium model based on the elastic wave theory, so as to understand the dispersion characteristics of the Rayleigh wave in the specific medium model. For the case in a layered medium, there exist a huge number of works in the literature. The classical propagator matrix method \cite{haskell1953} is the first systematical and efficient method for calculating the velocity of surface waves and the corresponding eigenvalue problems. However, this approach is unstable. That is, the result by this method is inaccurate in the case of high frequencies. In order to overcome this defect, many improved approaches, which are usually based on the propagator matrix method, had been developed. For example, the Schwab-Knopoff method \cite{knopoff1964}, the Abo-Zena method \cite{abo1979dispersion, menke1979comment}, etc. The idea of these methods is to improve the transfer matrix of the displacement stress vector, so as to eliminate the exponential term in the transfer process and obtain accurate results. Another line of research, initiated by \cite{kennett1979seismic} and named as the reflection and transmission coefficient method, has recently gained considerable attention. Compared with the transfer matrix based approaches, the reflection and transmission coefficient method is stable and very easily implemented in practice. However, when solving the dispersion equation in the reflection and transmission coefficient method, there exists the issue of ``mistake solution'' \cite{xia2006simple, xia2015new, cui2021velocity}. In order to solve this weakness, by combining the propagator matrix method and the reflection and transmission coefficient method, we proposed a new Rayleigh surface wave dispersion forward model in this work, which exhibits as a robust method for estimating dispersion characteristics such as  the traveling surface wave propagation constant, etc.

The inverse problem of estimating the medium parameters of Earth from the measured characteristic of waves (e.g., the dispersion curve) is an important issue in Geophysics \cite{xia2015new,cui2021velocity}.  Loosely speaking, numerical methods for solving such kind of inverse problems are usually classified into three groups: the linearization methods, the heuristic global-optimization methods, and the data-driven methods. The output of the linearization methods is very sensitive to the initial choice of the solution and the accuracy of the computed Jacobian matrix. Therefore, this type of methods often fails in the real world application. The heuristic global-optimization methods, such as the annealing simulation algorithm \cite{sun2017nonlinear},  the particle-swarm global optimization method \cite{poormirzaee2016s}, the genetic algorithm \cite{dal2007rayleigh,lei2019inversion,jianxun2018}, and others \cite{zhang2019wave,aleardi2020transdimensional,sivaram2018shear,mora2018inversion}, are usually derivative-free methods and the implementation is simple. Hence, during last decades, they have been widely used in the practice. The main defect of the heuristic global-optimization methods is their slow performance. Even for a small size problem, they usually need to take weeks to find the correct global minimizer.

The data-driven methods, especially the  machine learning methods, can serve as a remedy to all the above problems and begin to show potentials in geoscience. For example, \cite{LSD2019} propose a method that applies the Convolutional Neural Networks (CNN) and the principle component analysis to parameterize complex geological models in low dimensions. \cite{WangY2018} apply CNN to reconstruct high-resolution porous structures in rock samples from low-resolution $\mu$-CT images and high-resolution scanning electron microscope images. \cite{Holm-Jensen2020} approach a linear probabilistic waveform inversion problem in cross-hole tomography with machine learning techniques. In their study, the forward model is approximated with a linear operator using the ridge regression algorithm, which helps to find an analytical solution to the inversion problem. \cite{Mosser2020} use generative adversarial networks that represent prior distributions of the geological heterogeneities in their study on the stochastic seismic waveform inversion problem. \cite{HeWang2021} reparameterize model parameters in the full-waveform inversion problem as the weights of a deep neural network, which is known to have a strong representation capacity. For more recent progresses on deep learning methods in Geophysics, we refer to the latest review paper by \cite{YuMa2021} and references therein.

In the field of surface-wave tomography studies, applications of machine learning are still limited, yet they prove successful. For example, \cite{Meier2007} use the feed-forward neural network, and \cite{Cheng2019} use stacked auto-encoders, to invert surface wave data for crustal thicknesses. \cite{hu2020using} use convolutionary neural networks to invert the theoretical Rayleigh-wave phase and group velocity images for velocities. However, the models in these studies are incapable to deal with a significant problem that is prevalent in the inverse problems, which is called the non-uniqueness problem (see e.g., \cite{Shapiro2002}). That is, one dispersion curve may correspond to several distinct medium models, which leads to great difficulties in solving the inverse problem. Hence, the second goal in this work is to design a new neutral network architecture that can solve the non-uniqueness issue for our geophysics problem. The design of our new neutral network is based on the model of Mixture Density Network (MDN), which was firstly introduced by \cite{Bishop94,Bishop95} and now becomes widely used in various fields, including text generation (\cite{graves2014generating}), image recognition (\cite{Gkioxari2018}), unsupervised distribution and density estimations (\cite{Uria2018}), and transportation forecasting (\cite{Xu2018}).

The remainder of this paper is structured as follows. In Sect. \ref{sec:2}, we introduce a new recurrence equation for the dispersion of Rayleigh wave in a homogeneous and layered media. Moreover, in this section, we formulate the forward and inverse problems by a rigorous mathematical language. In  Sect. \ref{sec:3}, based on the mixture density network, we develop a new marching learning method for solving our inverse problem in geophysics. A detailed construction procedure and numerical experiments of our marching learning method are also reported in Sect. \ref{sec:3}. Finally, concluding remarks are given in Sect. \ref{Conclusion}.

\section{A novel recurrence equation for the dispersion of Rayleigh wave in a homogeneous and layered media}
\label{sec:2}

\begin{figure}[!b]
\includegraphics[scale=.6]{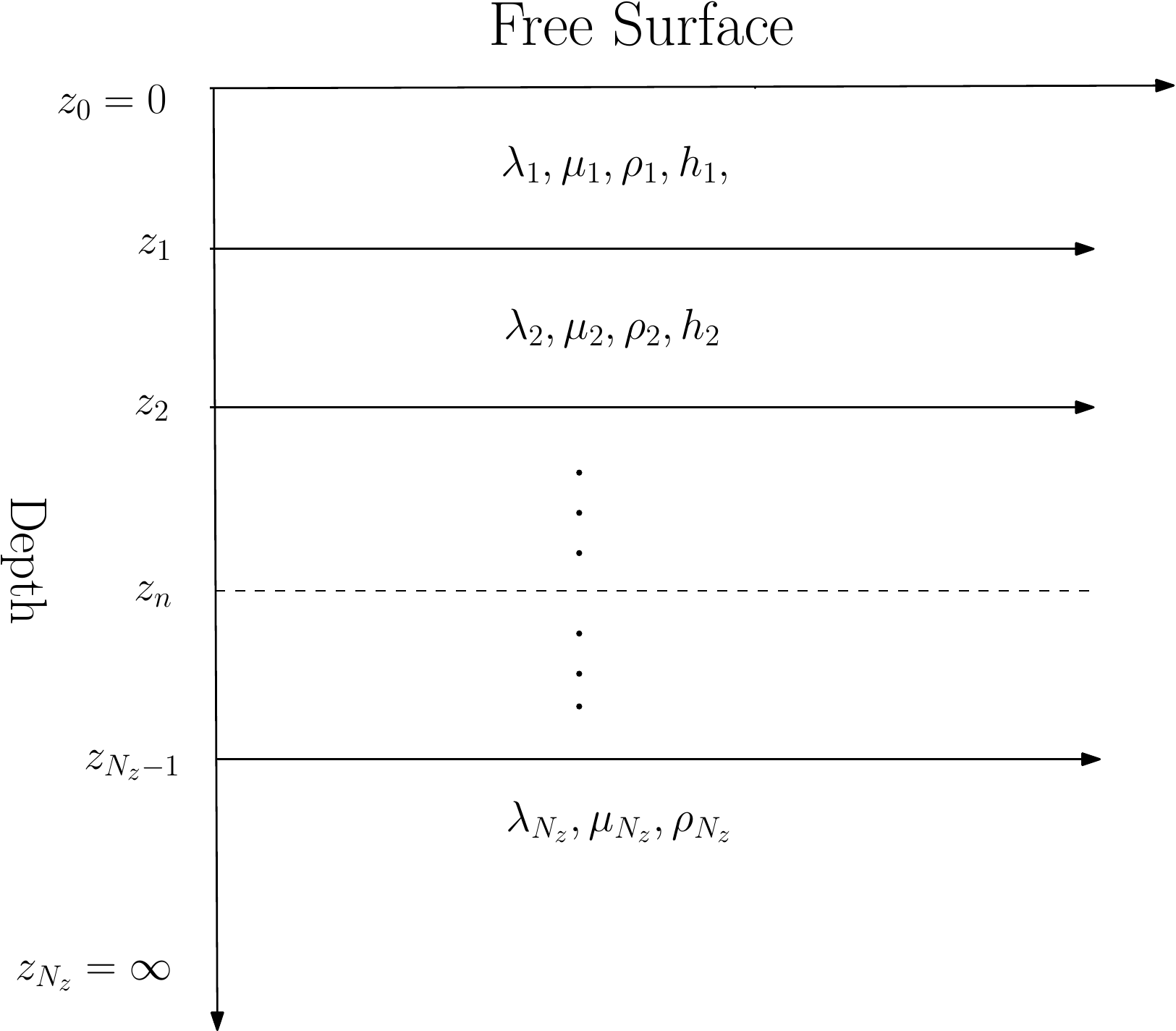}
\caption{The layered medium model of the Earth.}
\label{EarthModel}
\end{figure}

In this section, by introducing the quantity -- the seismic impedance tensor, we develop a new mathematical model for the efficient computation of the traveling wave characteristic in a layered medium. To that end, let us consider a layered medium with the boundary $z=z_n$, $n\in\{1,\cdots, N_z\}$, where $z_0=0$ represents the Earth's surface. In our simplified model (see Figure \ref{EarthModel}), inside each layer the Lam\'e coefficients $\lambda$ and $\mu$, as well as the Earth's density $\rho$ are all assumed to be constants. Specifically, for $z\in[z_{n-1},z_n]$, $\lambda\equiv \lambda_n$, $\mu\equiv \mu_n$ and $\rho\equiv \rho_n$. Moreover, the Cartesian coordinate system is constructed such that the traveling surface wave propagates along the $OX$ axis in this medium. Then, the wave displacement vector $\vec{U}$ has the form
\begin{equation}
 \label{example1}
 \vec{U}(x,y,z)=\overline{u}(z)e^{i \gamma x+i \omega t},
\end{equation}
where $i=\sqrt{-1}$ denotes the imaginary unit, $\omega$ denotes the frequency of the seismic field, $\gamma$ represents the traveling surface wave propagation constant (characteristic), and $\overline{u}(z)=(u_x,u_y,u_z)$ is the traveling-surface wave amplitude. By the assumption in our simplified model, inside a layer (e.g. $z\in[z_{n-1},z_n]$), $\vec{U}=(U_x,U_y,U_z)$ satisfies the Lam\'e equation with constant coefficients (\cite{Dmitriev})
\begin{equation}
 \label{example2}
(\lambda_n+2\mu_n){\rm grad} \, {\rm div}\vec{U}-\mu_n {\rm rot} \, {\rm rot} \vec{U}+\omega^2\rho_n\vec{U}=0, \quad z\in[z_{n-1},z_n].
\end{equation}

For Rayleigh waves, we only need to consider the $OX$ and $OZ$ components. By combining Eq. (\ref{example2}) and the stress equation (\cite{abo1979dispersion})
\begin{equation*}
\sigma_{ij}=\lambda \delta_{ij}{\rm div}\vec{U}+\mu \left(\frac{\partial U_i}{\partial j}+\frac{\partial U_j}{\partial i} \right), \quad {\rm for~all~} i,j\in\{x,y,z\},
\end{equation*}
we have
\begin{equation}
 \label{example3}
\begin{cases}
 \ (\lambda_n+2\mu_n)i\gamma (i\gamma u_x(z)+{u_z'})-\mu_n(-u_x''+i\gamma u_z')+\omega^2\rho_n u_x=0, \\
  (\lambda_n +2\mu_n)(i\gamma u_x'(z)+{u_z''}) - i\gamma\mu_n(u_x'+i\gamma u_z) + \omega^2\rho_n u_z=0, \\
 \ \sigma_{xz}=\mu_n(i\gamma u_z(z)+u_x'(x)), \\
 \ \sigma_{zz}=(\lambda_n +2\mu_n)u_z'(z)+i\gamma\lambda_n u_x(x),
  \end{cases}
\end{equation}
where $\delta_{ij}$ denotes the Kronecker delta function and $(\cdot)'=\frac{d}{dz}(\cdot)$ refers to the derivative with respect to the variable $z$. The boundary conditions for (\ref{example3}) in layered media are:
\begin{itemize}
 \item [1)]  On the earth's surface, the stress tensor is zero.
 \begin{equation}
 \label{example4}
  \sigma_{xz}|_{z=0}=\sigma_{zz}|_{z=0}=0.
\end{equation}
 \item [2)]Displacements and tensions are continuous for the j-th layer. That is,
 \begin{equation}
 \begin{aligned}
 \label{example5}
 u_x(z_j)|_{z=z_j+0}&= u_x(z_j)|_{z=z_j-0},&\\
 \ u_z(z_j)|_{z=z_j+0}&= u_z(z_j)|_{z=z_j-0},\\
 \ \sigma_{xz}(z_j)|_{z=z_j+0}&=\sigma_{xz}(z_j)|_{z=z_j-0},\\
 \ \sigma_{zz}(z_j)|_{z=z_j+0}&=\sigma_{zz}(z_j)|_{z=z_j-0}.\\
 \end{aligned}
\end{equation}
\item [3)] The vanishing boundary condition, which is
\begin{equation}
\label{vanishBoundary}
\lim_{z\to \infty} u_x(z) = \lim_{z\to \infty}u_z(z) =0.
\end{equation}
\end{itemize}

In this work, we assume that the seismic impedance tensor of a Rayleigh wave links the stresses $\sigma_{xz}$, $\sigma_{zz}$with the displacements $u_x$, $u_z$ by linear relationships
\begin{equation}
 \label{example6}
\begin{cases}
 \ \sigma_{xz}=Z_{xx} u_x + Z_{xz} u_z,\\
 \ \sigma_{zz}=Z_{zx} u_x + Z_{zz} u_z,
  \end{cases}
\end{equation}
where $Z_{xx}, Z_{xz}, Z_{zx}$ and $Z_{zz}$ are four constants, which should be determined later. An explanation of (\ref{example6}) can be found in the work of \cite{Dmitriev,jianxun2018,jianxun2019}. To this end, define the seismic impedance tensor $\hat Z$ of rank 2 by
\begin{equation}
 \label{example7}
\hat Z=\begin{pmatrix} Z_{xx} & Z_{xz} \\ Z_{zx} & Z_{zz} \end{pmatrix}.
\end{equation}

By Eq. (\ref{example3}) and the continuation condition (\ref{example5}) one can conclude that the impedance tensor $\hat Z$ is continuous with respect to parameters $\lambda$, $\mu$ and $\rho$.

Considering the impedance tensor on the surface of Earth, we derive  from the boundary condition (\ref{example4}) for the stress that (denote as $\hat Z_0= \hat Z (z = 0)$)
\begin{equation}
\begin{aligned}
 \label{example8}
\ 0= \sigma_{xz}(z=0)=\overset{(0)}{Z_{xx}} u_x(z=0)+\overset{(0)}{{Z_xz}} u_z(z=0),\\
 \ 0= \sigma_{zz}(z=0)=\overset{(0)}{Z_{zx}} u_x(z=0)+\overset{(0)}{{Z_zz}} u_z(z=0).
\end{aligned}
\end{equation}
Hence, a Rayleigh wave exists (i.e., $u_x (0)u_z (0) \neq 0$) if and only if the following equation holds
\begin{equation}
 \label{example9}
\det \hat Z_0=\overset{(0)}{Z_{xx}}\overset{(0)}{Z_{zz}}-\overset{(0)}{Z_{xz}}\overset{(0)}{Z_{zx}}=0.
\end{equation}

The essential idea in this work is to use the dispersion equation (\ref{example9}) to establish the new relationship between the Rayleigh wave propagation constant $\gamma$ and other parameters such as the frequency $\omega$ and the parameters of the layered medium.

Denote by  $u_s (z)$ and $u_p (z)$ the solenoidal and potential parts of the displacement field, respectively (\cite{abo1979dispersion, haskell1953, knopoff1964}). Then, for a Rayleigh wave, we have
\begin{equation}
\begin{aligned}
 \label{example10}
u_x&=-u_s'(z)+i\gamma u_p(z),\\
u_z&=i\gamma u_s(z)+u_p'(z),
\end{aligned}
\end{equation}
where $u_s (z)$ is the solution of equation $u_s''-\eta_s^2 u_s=0$. Here, $\eta_s=\sqrt{\gamma^2-k_s^2}$, $k_s=\frac{\omega}{v_s}$, $\rm Re~ \eta_s>0$, $v_s=\sqrt{\mu / \rho}$ is the transverse wave velocity. $u_p (z)$ is the solution of equation $u_p''-\eta_p^2 u_p=0$, where $\eta_p=\sqrt{\gamma^2-k_p^2}$, $k_p=\frac{\omega}{v_p}$, $\rm Re~ \eta_p>0$, $v_p=\sqrt{(\lambda+2\mu) / \rho}$ is the longitudinal wave velocity. Note that in both cases $\gamma >k_s>k_p$.

Consider the field in layer $z\in[z_{n-1},z_n]$ with the constants $\lambda_n$, $\mu_n$, $\rho$, where $n\in\left\{ 1, \cdots, N_z \right\}$, $z_0=0$, $z_{N_z}=-\infty$, and $h_n=z_n-z_{n-1}$ ($n\in\left\{ 1, \cdots, N_z-1 \right\}$) is the layer thickness. According to the general solution of the equation $u''-\eta^2 u=0$, in layer $n$, $u_s(z)$ and $u_p(z)$ can be represented in the form
\begin{equation}
\begin{aligned}
 \label{example11}
u_s(z)=\overset{(n)}{A_s}e^{\overset{(n)}{\eta_s}(z_n-z_{n-1})}+\overset{(n)}{B_s}e^{-\overset{(n)}{\eta_s}(z_n-z_{n-1})},\\
u_p(z)=\overset{(n)}{A_p}e^{\overset{(n)}{\eta_p}(z_n-z_{n-1})}+\overset{(n)}{B_p}e^{-\overset{(n)}{\eta_p}(z_n-z_{n-1})},
\end{aligned}
\end{equation}
where $\overset{(n)}{A_s}, \overset{(n)}{A_p}, \overset{(n)}{B_s}$ and $\overset{(n)}{B_p}$ are four constants that depend on the boundary condition and will be determined later.

It follows from (\ref{example10}) and (\ref{example11}) that in the $n$-th layer, we can get the expression of displacement and stress tensor
\begin{equation}
 \label{example12}
\begin{cases}
 \ u_x(z)=i\gamma\left( \overset{(n)}{A_p}e^{\overset{(n)}{\eta_p}(z-z_{n})}+\overset{(n)}{B_p}e^{-\overset{(n)}{\eta_p}(z-z_{n})}\right)
 -\overset{(n)}{\eta_s}\left(\overset{(n)}{A_s}e^{\overset{(n)}{\eta_s}(z-z_{n})}-\overset{(n)}{B_s}e^{-\overset{(n)}{\eta_s}(z-z_{n})}\right),\\
  \ u_z(z)=i\gamma\left( \overset{(n)}{A_s}e^{\overset{(n)}{\eta_s}(z-z_{n})}+\overset{(n)}{B_s}e^{-\overset{(n)}{\eta_s}(z-z_{n})}\right)
 +\overset{(n)}{\eta_p}\left(\overset{(n)}{A_p}e^{\overset{(n)}{\eta_p}(z-z_{n})}-\overset{(n)}{B_p}e^{-\overset{(n)}{\eta_p}(z-z_{n})}\right),\\
 \
 \sigma_{xz}(z)=\mu_n \Bigg(
 -\left(\gamma^2-\overset{(n)}{\eta_s}^2
 \right)\left(\overset{(n)}{A_s}e^{\overset{(n)}{\eta_s}(z-z_{n})}-\overset{(n)}{B_s}e^{-\overset{(n)}{\eta_s}(z-z_{n})}\right)\\
  \qquad\qquad  +2i\gamma\overset{(n)}{\eta_p}\left(\overset{(n)}{A_p}e^{\overset{(n)}{\eta_p}(z-z_{n})}-\overset{(n)}{B_p}e^{-\overset{(n)}{\eta_p}(z-z_{n})} \right)
 \Bigg),
\\
 \ \sigma_{zz}(z)=2i\gamma\mu_n\overset{(n)}{\eta_s}\left(\overset{(n)}{A_s}e^{\overset{(n)}{\eta_s}(z-z_{n})}-\overset{(n)}{B_s}e^{-\overset{(n)}{\eta_s}(z-z_{n})}\right)\\
   \qquad\qquad +\left((\lambda_n+2\mu_n)\overset{(n)}{\eta_p}^2-\lambda_n\gamma^2\right)\left(\overset{(n)}{A_p}e^{\overset{(n)}{\eta_p}(z-z_{n})}-\overset{(n)}{B_p}e^{-\overset{(n)}{\eta_p}(z-z_{n})}\right).
  \end{cases}
\end{equation}

According to the upper and lower interfaces of $n$-th layers in (\ref{example5}), we can derive together with the impedance tensor (\ref{example6}) and Eq. (\ref{example12}) that
\begin{equation}
 \label{example13}
\begin{cases}
 \ u_x(z_{n-1})-i\gamma\left( \overset{(n)}{A_p}e^{\overset{(n)}{\eta_p}(-h)}+\overset{(n)}{B_p}e^{-\overset{(n)}{\eta_p}(-h)}\right)
 +\overset{(n)}{\eta_s}\left(\overset{(n)}{A_s}e^{\overset{(n)}{\eta_s}(-h)}-\overset{(n)}{B_s}e^{-\overset{(n)}{\eta_s}(-h)}\right)=0\\
  \ u_z(z_{n-1})-i\gamma\left( \overset{(n)}{A_s}e^{\overset{(n)}{\eta_s}(-h)}+\overset{(n)}{B_s}e^{-\overset{(n)}{\eta_s}(-h)}\right)
 -\overset{(n)}{\eta_p}\left(\overset{(n)}{A_p}e^{\overset{(n)}{\eta_p}(-h)}-\overset{(n)}{B_p}e^{-\overset{(n)}{\eta_p}(-h)}\right)=0\\
 \
 \mu_n \Bigg(
 -\left(\gamma^2-\overset{(n)}{\eta_s}^2
 \right)\left(\overset{(n)}{A_s}-\overset{(n)}{B_s}\right)+2i\gamma\overset{(n)}{\eta_p}\left(\overset{(n)}{A_p} -\overset{(n)}{B_p} \right)n\Bigg) \\
 \qquad\qquad\qquad\qquad\qquad\qquad\qquad\qquad\qquad - \overset{(n)}{Z_{xx}}(z_n)u_x(z_n)-\overset{(n)}{Z_{xz}}(z_n)u_z(z_n)=0
\\

 \ 2i\gamma\mu_n\overset{(n)}{\eta_s}\left(\overset{(n)}{A_s}-\overset{(n)}{B_s}\right)
  +\left((\lambda_n+2\mu_n)\overset{(n)}{\eta_p}^2-\lambda_n\gamma^2\right)\left(\overset{(n)}{A_p}-\overset{(n)}{B_p}\right) \\
 \qquad\qquad\qquad\qquad\qquad\qquad\qquad\qquad\qquad
  -\overset{(n)}{Z_{zx}}(z_n)u_x(z_n)-\overset{(n)}{Z_{zz}}(z_n)u_z(z_n)=0
  \end{cases}
\end{equation}

For convenience, we rewrite the above system of equations in the form of matrix equation
\begin{equation}
 \label{example14}
G_n
\begin{pmatrix}
\overset{(n)}{A_s}\\ \overset{(n)}{B_s}\\ \overset{(n)}{A_p}\\ \overset{(n)}{B_p}
\end{pmatrix}=
\begin{pmatrix}
\overset{(n)}{G_{11}}&\overset{(n)}{G_{12}}&\overset{(n)}{G_{13}}&\overset{(n)}{G_{14}}\\ \overset{(n)}{G_{21}}&\overset{(n)}{G_{22}}&\overset{(n)}{G_{23}}&\overset{(n)}{G_{24}}\\ \overset{(n)}{G_{31}}&\overset{(n)}{G_{32}}&\overset{(n)}{G_{33}}&\overset{(n)}{G_{34}}\\ \overset{(n)}{G_{41}}&\overset{(n)}{G_{42}}&\overset{(n)}{G_{43}}&\overset{(n)}{G_{44}}
\end{pmatrix}
\begin{pmatrix}
\overset{(n)}{A_s}\\ \overset{(n)}{B_s}\\ \overset{(n)}{A_p}\\ \overset{(n)}{B_p}
\end{pmatrix}=
\begin{pmatrix}
u_x(z_{n-1})\\ u_z(z_{n-1})\\0\\ 0
\end{pmatrix},
 \end{equation}
 where
 \begin{equation}
 \label{example15}
 \begin{cases}
\overset{(n)}{G_{11}}=\overset{(n)}{\eta_s}e^{-\overset{(n)}{\eta_s}h}, ~
\overset{(n)}{G_{12}}=-\overset{(n)}{\eta_s}e^{\overset{(n)}{\eta_s}h}, ~
\overset{(n)}{G_{13}}=-i\gamma e^{-\overset{(n)}{\eta_p}h}, ~
\overset{(n)}{G_{14}} =-i\gamma e^{\overset{(n)}{\eta_p}h}, \\
\overset{(n)}{G_{21}}=-i\gamma e^{\overset{(n)}{\eta_s}h}, ~
\overset{(n)}{G_{22}}=-i\gamma e^{\overset{(n)}{\eta_s}h}, ~
\overset{(n)}{G_{23}}=-\overset{(n)}{\eta_p}e^{-\overset{(n)}{\eta_p}h}, ~
\overset{(n)}{G_{24}}=\overset{(n)}{\eta_p}e^{-\overset{(n)}{\eta_p}h}, \\
\overset{(n)}{G_{31}}=-\overset{(n)}{Z_{xx}}\overset{(n)}{\eta_s}+i\overset{(n)}{Z_{xz}}\gamma-\mu_n(-\overset{(n)}{\eta_s}^2-\gamma^2), ~
\overset{(n)}{G_{32}}=\overset{(n)}{Z_{xx}}\overset{(n)}{\eta_s}+i\overset{(n)}{Z_{xz}}\gamma-\mu_n(-\overset{(n)}{\eta_s}^2-\gamma^2)\\
\overset{(n)}{G_{33}}=-2i\mu_n\gamma\overset{(n)}{\eta_p}+i\overset{(n)}{Z_{xx}}\gamma+\overset{(n)}{Z_{xz}}\overset{(n)}{\eta_p}, ~
\overset{(n)}{G_{34}}=2i\mu_n\gamma\overset{(n)}{\eta_p}+i\overset{(n)}{Z_{xx}}\gamma-\overset{(n)}{Z_{xz}}\overset{(n)}{\eta_p}\\
\overset{(n)}{G_{41}}=-\overset{(n)}{Z_{zx}}\overset{(n)}{\eta_s}+i\overset{(n)}{Z_{zz}}\gamma-i(\lambda_n+2\mu_n)\gamma\overset{(n)}{\eta_s}+i\lambda_n\gamma\overset{(n)}{\eta_s}, \\
\overset{(n)}{G_{42}}=\overset{(n)}{Z_{zx}}\overset{(n)}{\eta_s}+i\overset{(n)}{Z_{zz}}\gamma+i(\lambda_n+2\mu_n)\gamma\overset{(n)}{\eta_s}-i\lambda_n\gamma\overset{(n)}{\eta_s},\\
\overset{(n)}{G_{43}}=i\overset{(n)}{Z_{zx}}\gamma+\overset{(n)}{Z_{zz}}\overset{(n)}{\eta_p}-(\lambda_n+2\mu_n)\overset{(n)}{\eta_p}^2+\lambda_n\gamma^2,\\
\overset{(n)}{G_{44}}=i\overset{(n)}{Z_{zx}}\gamma-\overset{(n)}{Z_{zz}}\overset{(n)}{\eta_p}-(\lambda_n+2\mu_n)\overset{(n)}{\eta_p}^2+\lambda_n\gamma^2.
  \end{cases}
 \end{equation}

From the system of linear algebraic equations (\ref{example14}), we can formally obtain
\begin{equation}
 \label{example16}
\begin{matrix}
\begin{pmatrix}
\overset{(n)}{A_s}\\ \overset{(n)}{B_s}\\ \overset{(n)}{A_p}\\ \overset{(n)}{B_p}
\end{pmatrix}
=G^{-1}_n
\begin{pmatrix}
u_x(z_{n-1})\\ u_z(z_{n-1})\\0\\ 0
\end{pmatrix}.
\end{matrix}
\end{equation}

On the other hand, we rewrite the last two equations of (\ref{example12}) in the following vector form:
\begin{equation}
 \label{example17}
 \begin{pmatrix}
 \sigma_{xz}(z_{n-1})\\\sigma_{zz}(z_{n-1})
 \end{pmatrix}
 =\Gamma_n \cdot\rm diag\begin{pmatrix}
 e^{\overset{(n)}{\eta_s}h},&e^{-\overset{(n)}{\eta_s}h},&e^{\overset{(n)}{\eta_p}h},&e^{-\overset{(n)}{\eta_p}h}
 \end{pmatrix}
 \begin{pmatrix}
 \overset{(n)}{A_s}\\ \overset{(n)}{B_s}\\ \overset{(n)}{A_p}\\ \overset{(n)}{B_p}
 \end{pmatrix},
 \end{equation}
where

\begin{equation*}
 \label{example18}
\Gamma_n =
 \begin{pmatrix}
 \mu_n(-(\overset{(n)}{\eta_s})^2-\gamma^2)&
 \mu_n(-(\overset{(n)}{\eta_s})^2-\gamma^2)&
 2i\mu_n\gamma\overset{(n)}{\eta_p}&
 -2i\mu_n\gamma\overset{(n)}{\eta_p}\\
2i\mu_n\gamma\overset{(n)}{\eta_s}&
2i\mu_n\gamma\overset{(n)}{\eta_s}&
(\lambda_n+2\mu_n)(\overset{(n)}{\eta_p})^2-\lambda_n\gamma^2&
(\lambda_n+2\mu_n)(\overset{(n)}{\eta_p})^2-\lambda_n\gamma^2
 \end{pmatrix} .
\end{equation*}

Substituting (\ref{example16}) into (\ref{example17}), we get
\begin{equation}
 \label{example19}
\begin{pmatrix}
 \sigma_{xz}(z_{n-1})\\\sigma_{zz}(z_{n-1})
 \end{pmatrix}
 =\Gamma_n \cdot\rm diag\begin{pmatrix}
 e^{\overset{(n)}{\eta_s}h},&e^{-\overset{(n)}{\eta_s}h},&e^{\overset{(n)}{\eta_p}h},&e^{-\overset{(n)}{\eta_p}h}
 \end{pmatrix}
G^{-1}_n
\begin{pmatrix}
u_x(z_{n-1})\\ u_z(z_{n-1})\\0\\ 0
\end{pmatrix}.
 \end{equation}
Define
\begin{equation}
\label{Edecomposition}
( E_n | \hat{E}_n) := \Gamma_n \cdot\rm diag\begin{pmatrix}
 e^{\overset{(n)}{\eta_s}h},&e^{-\overset{(n)}{\eta_s}h},&e^{\overset{(n)}{\eta_p}h},&e^{-\overset{(n)}{\eta_p}h}
 \end{pmatrix}
G^{-1}_n,
\end{equation}
where $E_n$ and $\hat{E}_n$ are both matrices of size $2\times2$. By the definition of $E_n$, it is a nonlinear function of parameters $\overset{(n)}{Z_{xx}}, \overset{(n)}{Z_{xz}}, \overset{(n)}{Z_{zx}}$ and $\overset{(n)}{Z_{zz}}$ (or equivalently the matrix $\hat{Z}$ defined in (\ref{example7})) at the $n$-th layer. That is, 
\begin{equation}
 \label{ZxzNew}
E_n = E_n(\hat{Z}_{n}), \qquad
\hat{Z}_{n} =\begin{pmatrix} \overset{(n)}{Z_{xx}} & \overset{(n)}{Z_{xz}} \\ \overset{(n)}{Z_{zx}} &\overset{(n)}{ Z_{zz}} \end{pmatrix}.
\end{equation}

Furthermore, by using (\ref{example19}), (\ref{Edecomposition}) and (\ref{ZxzNew}) we deduce that
\begin{equation}
 \label{NewEq1}
\begin{pmatrix}
 \sigma_{xz}(z_{n-1})\\\sigma_{zz}(z_{n-1})
 \end{pmatrix}
 = E_n(\hat{Z}_{n})
\begin{pmatrix}
u_x(z_{n-1})\\ u_z(z_{n-1})
\end{pmatrix}.
 \end{equation}

Form the relation (\ref{example6}), we know that the impedance tensor at the $(n-1)$-th layer has the representation
\begin{equation}
 \label{NewEq2}
\begin{pmatrix}
 \sigma_{xz}(z_{n-1})\\\sigma_{zz}(z_{n-1})
 \end{pmatrix}
 = \hat{Z}_{n-1}
\begin{pmatrix}
u_x(z_{n-1})\\ u_z(z_{n-1})
\end{pmatrix}.
 \end{equation}

Combine (\ref{NewEq1}) and (\ref{NewEq2}) to obtain that $\hat{Z}_{n-1} = E_n(\hat{Z}_{n})$, which implies that
\begin{equation}
 \label{EnZn}
\hat{Z}_{0} = E_1(\hat{Z}_{1}) = E_1(E_2(\hat{Z}_{2})) = \cdots =  E_1(E_2( \cdots E_{N_z}(\hat{Z}_{N_z}))).
\end{equation}
Thus, Eq. (\ref{example9}) gives
\begin{equation}
 \label{detZ0}
\det E_1(E_2( \cdots E_{N_z}(\hat{Z}_{N_z}))) = 0.
\end{equation}

Note that the quantity $\hat{Z}_{N_z}$ is a matrix-valued function of parameters $\gamma$ and $\mu_{N_z}$, $\lambda_{N_z}$, $\rho_{N_z}$. Actually, it follows\footnote{It can be obtained by using (\ref{example3}), (\ref{example6}), and the vanishing boundary conditions (\ref{vanishBoundary})} from \cite[Formula (42) ]{jianxun2018} that
\begin{equation}
\begin{aligned}
 \label{EnZn2}
& \hat{Z}_{N_z} = \hat{Z}_{N_z} (\gamma)= \\ &~
\begin{pmatrix} -\frac{{\mu}_{N_z}\overset{({N_z})}{\eta_p}((\overset{({N_z})}{\eta_s})^2-\gamma^2)}{\overset{({N_z})}{\eta_s}\overset{({N_z})}{\eta_p}-\gamma^2} &\qquad&
\frac{i{\mu}_{N_z}\gamma((\overset{({N_z})}{\eta_s})^2-2\overset{({N_z})}{\eta_s}\overset{({N_z})}{\eta_p}+\gamma^2))}{\overset{({N_z})}{\eta_s}\overset{({N_z})}{\eta_p}-\gamma^2} \\
&&\\
({\lambda}_{N_z}+2{\mu}_{N_z})\frac{i\gamma(\overset{({N_z})}{\eta_s}\overset{({N_z})}{\eta_p}-(\overset{({N_z})}{\eta_p})^2)}{\overset{({N_z})}{\eta_s}\overset{({N_z})}{\eta_p}-\gamma^2}+i\gamma{\lambda}_{N_z} &\qquad&
 ({\lambda}_{N_z}+2{\mu}_{N_z})\frac{\overset{({N_z})}{\eta_s}((\overset{({N_z})}{\eta_p})^2-\gamma^2)}{\overset{({N_z})}{\eta_s}\overset{({N_z})}{\eta_p}-\gamma^2}
\end{pmatrix}.
\end{aligned}
\end{equation}

For simplicity, we assume that $h_n \equiv h\approx 4$km for all $n=1, \cdots, {N_z}-1$. Now, let us discuss the values of parameters $\mu_n, \lambda_n$ and $\rho_n$ in the $n$-th layer. To that end, denote by $V^{(n)}_S$ and $V^{(n)}_P$ the S-wave velocity and P-wave velocity in the $n$-th layer, respectively. By the the empirical formula (cf. \cite{zhu1995}) we have $\rho_n=0.466 [V^{(n)}_S]^{0.214}$. Moreover, it is well known that (see e.g., \cite{abo1979dispersion}), $V^{(n)}_S=\sqrt{\mu_n/\rho_n}$ and $V^{(n)}_P=\sqrt{(\lambda_n+2\mu_n)/\rho_n}$, which implies together with the empirical formula $V^{(n)}_P = 1.732V^{(n)}_S$ (see e.g. \cite{chai2019correlation}) and the mentioned relation $\rho_n=0.466[V^{(n)}_S]^{0.214}$ that $\lambda_n = \mu_n= 0.466 [V^{(n)}_S]^{2.214}$. To summarize, in this work, all parameters $\mu_n, \lambda_n$ and $\rho_n$ are purely dependent on the value of the S-wave velocity $V^{(n)}_S$.

In practice, the phase velocity $C$ of the Rayleigh wave at a given frequency $\omega$ is usually measured. it is well known that $C=\omega/\gamma$ (\cite{abo1979dispersion}). To that end, denote by $\bar{V}_S = (V^{(1)}_S, \cdots , V^{({N_z})}_S)$ the vector of S-wave velocity in our layered medium model of the Earth. By the construction of matrix $\{E_n\}^{N_z}_{n=1}$, we conclude that the left-hand side in (\ref{detZ0}) is a nonlinear continuous function of parameters $C$ and $\bar{V}_S$. According to the implicit function theorem, we deduce from Eq. (\ref{detZ0}) that in some neighborhood of the point $\bar{V}_S$, there exists a differentiable function $f$ of $\bar{V}_S$ such that
\begin{equation}
 \label{f}
f(\bar{V}_S) = C.
\end{equation}

Hence, for a given S-wave velocity $\bar{V}_S$, one can efficiently compute the phase velocity $C$ through the procedure (\ref{f}). It should be noted that although we do not obtain the closed form of function $f$, the value of $C$ can be numerically determined by solving a nonlinear equation (\ref{detZ0}). Since both the physical quantity $C$ and the function $f$ depend on the frequency, given a frequencies $\{\omega_j\}^{N_{\omega}}_{j=1}$, one can obtain a sequence of equations
\begin{equation*}
 \label{fj}
f(\bar{V}_S;\omega_j) = C(\omega_j), \quad j=1, \cdots, N_{\omega}.
\end{equation*}
Define $\mathbf{f}=(f(\bar{V}_S;\omega_1), \cdots, f(\bar{V}_S;\omega_{N_\omega}))$ and $\mathbf{C}=(C(\omega_1), \cdots, C(\omega_{N_\omega}))$. Without ambiguity, we can call $\mathbf{C}$ the dispersion curve (actually, $\mathbf{C}$ represents the discretized version of the dispersion curve). Now, we obtain the vector form of the equation
\begin{equation}
 \label{fjVector}
\mathbf{f}(\bar{V}_S) = \mathbf{C}.
\end{equation}
In the next work, we will develop a robust machine learning method to estimate the S-wave velocity $\bar{V}_S$ at a given frequency from the dispersion curve and our new mathematical model (\ref{fjVector}).

\section{An efficient machine learning based inverse solver}
\label{sec:3}

This section is devoted to the construction of a neural network that serves as an efficient numerical solver for the inverse problem (\ref{fjVector}). In this solver, we assume that both of $\bar{V}_S\in \mathbb{R}^{N_z}$ and $\mathbf{C}\in \mathbb{R}^{N_{\omega}}$ are random vectors, where $N_z$ and $N_{\omega}$ are the total number of crustal layers and frequencies, respectively. The problem of estimating S-wave velocity $\bar{V}_S$ from the measured dispersion curve $\mathbf{C}$ is transferred to estimating the function from $\mathbf{C}$ to $\bar{V}_S$, which will be solved by using the models based on neural networks in this work.

\subsection{Data structure}
We use neural networks to approach each of the 3-layer, 5-layer and 9-layer case (i.e., $N_z=3,5,9$) separately. Let $(\mathbf{x}_j, \mathbf{y}_j):=(\{x_{lj}\}_{l=0}^{N_z-1}, \{y_{lj}\}_{l=0}^{N_{\omega}-1})$ represent the $j^{th}$ realization of $(\bar{V}_S,\mathbf{C})$ in the $N_z$-layer case, where we take $N_{\omega}=50$ frequencies uniformly distributed in $[0.0785, 12.57]$, and the $N$ samples\footnote{We generate $48,000$ for each of the 3-layer and 5-layer case, and generate $120,000$ models for the 9-layer case as it is more complicated.} are generated in the following way.
\begin{itemize}
\item For each of $l\in\{0,1,\cdots,N_z-1\}$ and $j\in\{0,1,\cdots,N-1\}$, $x_{lj}$ is sampled independently from Uniform($a_l$,$b_l$), where $a_l$ and $b_l$ are taken from the global crustal average model (\cite{laske2013update}) and their values are listed in Table \ref{scope}.

\begin{table}[!htb]
\caption{Ranges for randomly generating entries in $\mathbf{x}$.}
\begin{center}
\begin{tabular}{p{2cm} p{2.5cm}  p{2.5cm} p{2.5cm} }
\toprule
                        &3-layer             &5-layer                & 9-layer \\
\toprule
$[a_0, b_0]$    &$[3.00,4.00]$       &$[3.00,3.80]$          &$[3.00,3.80]$   \\
$[a_1, b_1]$    &$[3.80,4.80]$       &$[3.20,4.00]$          &$[3.10,3.90]$   \\
$[a_2, b_2]$    &$[4.60,5.60]$       &$[3.80,4.60]$          &$[3.20,3.95]$   \\
$[a_3, b_3]$    &  -                     &$[3.80,4.60]$          &$[3.30,4.00]$   \\
$[a_4, b_4]$    &  -                     &$[4.00,4.80]$          &$[3.80,4.60]$   \\
$[a_5, b_5]$    &  -                     &-                          &$[3.90,4.70]$   \\
$[a_6, b_6]$    &  -                     &-                          &$[4.00,4.75]$   \\
$[a_7, b_7]$    &  -                     &-                          &$[4.20,4.80]$   \\
$[a_8, b_8]$    &  -                     &-                          &$[4.60,5.60]$   \\
\bottomrule
\end{tabular}
\label{scope}
\end{center}
\end{table}

\item The 50 entries of $\mathbf{y}_j$ are computed using the forward model (\ref{fjVector}) with the generated $\mathbf{x}_j$.
\end{itemize}

Moreover, in the simulation, we then randomly assign $80\%$ of the $N$ samples for model training, $10\%$ for validation, and $10\%$ for testing.

\subsection{Approach the inversion problem with a machine learning method based on Mixture Density Network (MDN)}

As demonstrated in Appendix A, due to the non-uniqueness for the inversion mapping of (\ref{fjVector}), the ordinary Feedforward Neural Networks are not suitable for the considered inverse problem. In this section, we will fix this problem by using the Mixture Density Network (MDN).

\subsubsection{A brief introduction to MDN}
\label{introMDN}
Suppose that the dataset $\{(\mathbf{x}_j,\mathbf{y}_j)\}_{j=0}^{N-1}$ is generated by an unknown multivalued function $f^{-1}(\mathbf{y})$ with $K$ branches $\{f^{-1}_k(\mathbf{y})\}_{k=0}^{K-1}$. Then it is a natural idea to characterize the data with a probability model. For example, for each sample $(\mathbf{x},\mathbf{y})$ it can be assumed that
\[ \mathbf{x} = \sum_{k=0}^{K-1} \pi_k f^{-1}_k(\mathbf{y}),\]
where $\pi_k\in[0,1]$ is the probability for $(\mathbf{x},\mathbf{y})$ to lie on the $k$-th branch.

The mixture model is based on a similar idea, but rather than determistic values it allows $\mathbf{x}$ on the $k$-th branch (for $k\in\{0,1,2,\cdots,K-1\}$) to be the realization of a certain distribution, say a Gaussian distribution with mean $f^{-1}_k(\mathbf{y})$. Specifically, for any sample $(\mathbf{x}, \mathbf{y})$, the Gaussian mixture (GM) model assumes that $\mathbf{x}$ is drawn from a conditional distribution $\mathbfit{X}|\mathbfit{Y}=\mathbf{y}$ with the following probability density function (pdf):
\begin{gather}
p_{\mathbfit{X} | \mathbfit{Y}=\mathbf{y}}(\mathbf{x}) = \sum_{k=0}^{K-1} \pi_k p_k(\mathbf{x})  \text{ for } k=0, 1, \cdots, K-1,
\end{gather}
where $p_k$ denotes the pdf of a Gaussian distribution $\mathcal{N}(\mu_k,\sigma_k^2)$, and the weight $\pi_k\in [0,1]$ satisfies $\sum_{k=0}^{K-1} \pi_k = 1$. Since all the parameters' values depend on $\mathbf{y}$, we can also write them as $\pi_k(\mathbf{y}), \mu_k(\mathbf{y}), \sigma_k^2(\mathbf{y})$. Note that the $k^{th}$ branch $f_k^{-1}(\mathbf{y})$ in the above deterministic case corresponds to $\mathcal{N}(\mu_k(\mathbf{y}),\sigma_k^2(\mathbf{y}))$. 

The Mixture Density Network (MDN), when the prior is chosen to be Gaussian, is a GM with one more assumption: $\pi_k(\mathbf{y}), \mu_k(\mathbf{y})$, and $\sigma_k^2(\mathbf{y})$ are continuous functions and hence they can be approximated with a feed-forward neural network. This model assumes that the data samples are realizations of two random vectors $\mathbfit{X}$ and $\mathbfit{Y}$, whose distributions are as follows:  conditional on $\mathbfit{Y}=\mathbf{y}$ the distribution of $\mathbfit{X}$ is a weighted mixture of $K$ Gaussian distributions. That is, the pdf of $\mathbfit{X}$ conditional on $\mathbfit{Y} = \mathbf{y}$ is
\[ p_{\mathbfit{X} | \mathbfit{Y}=\mathbf{y}}(\boldsymbol{x}) = \sum_{k=0}^{K-1} \pi_k p_k(\boldsymbol{x}), \]
where $p_k$ denotes the probability density function of $\mathcal{N}(\mu_k, \sigma_k^2)$ for $k=0,1,\cdots, K-1$, and the weight $\pi_k\in [0,1]$ satisfies $\sum_{k=0}^{K-1} \pi_k = 1$. The parameters, $\pi_k(\mathbf{y})$, $\mu_k(\mathbf{y})$, and $\sigma_k^2(\mathbf{y})$, are continuous functions on $\mathbf{y}$, and hence can be approximated with an FNN with input $\mathbf{y}$.

\subsubsection{MDN model structure}
\label{mdnstruct}
\begin{figure}[!htbp]
	\centering
		\includegraphics[scale=.6]{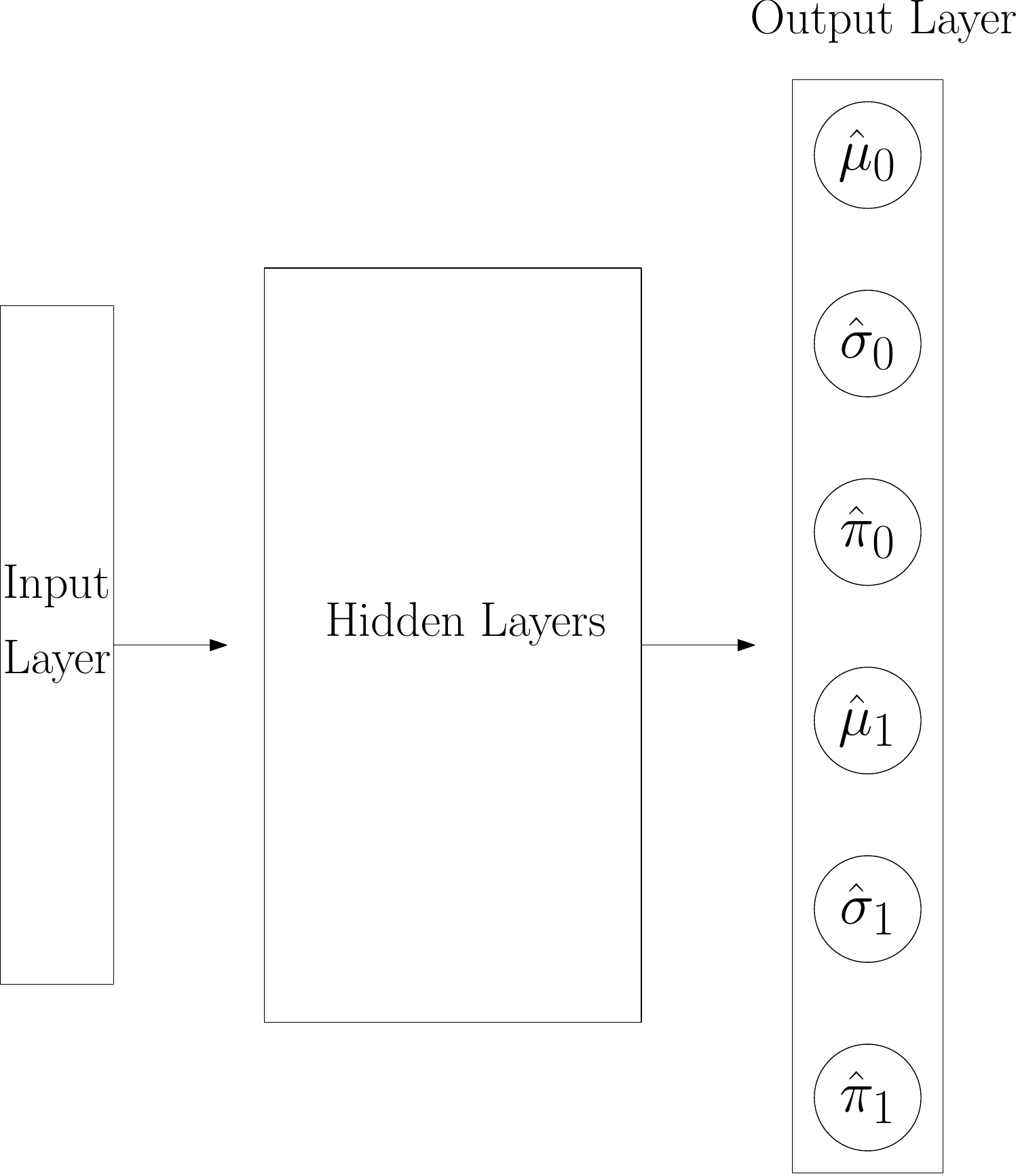}
       	\caption{Illustration of the FNN included in an MDN}
\label{MDN}
\end{figure}

The structure of an MDN is the same as an FNN whose input is $\mathbf{y}$ and whose target output is $\{\pi_k(\mathbf{y}), \mu_k(\mathbf{y}), \sigma_k^2(\mathbf{y})\}_{k=0}^{K-1}$. So the MDN output can be viewed as a Gaussian mixture distribution. Figure \ref{MDN} illustrates such an FNN with the simple case where the dimension of $\mathbf{y}$ is one and $K=2$.

For a more general case, suppose the dimension of $\mathbf{x}$ is $n$, then $\mu_k$ and $\sigma_k$ also share the same dimension. So the dimension of $(\mu_k, \sigma_k, \pi_k)$ is $2n+1$. Hence, the number of nodes in the output layer of the FNN is $(2n+1)K$. We index nodes in the output layer with $\{0,1,\cdots,(2n+1)K-1\}$ and divide them into K sets, and set the target of the $k$-th set as $(\pi_k(\mathbf{x}), \mu_k(\mathbf{x}), \sigma_k^2(\mathbf{x}))$. More specifically, let $\{o_j\}_{j=0}^{(2n+1)K-1}$ denote the FNN output for an input $\mathbf{y}$, then
\begin{itemize}
\item the target value for $\hat{\mu}_k(\mathbf{y}):=[o_{(2n+1)k}, o_{(2n+1)k+1}, \cdots, o_{(2n+1)k+(n-1)}]$ is $\mu_k(\mathbf{y})$, where a ReLU activation is used for these nodes to make sure the outputs are non-negative,
\item the target value for $\hat{\sigma}_k(\mathbf{y}):=[o_{(2n+1)k+n}, o_{(2n+1)k+n+1}, \cdots, o_{(2n+1)k+(2n-1)}]$ is $\sigma_k(\mathbf{y})$, where a sigmoid activation multiplied by a positive hyper-parameter\footnote{This hyper-paramter is set to adapt to variances of the data samples and its value is determined by the validation-set approach. We find in our later experiments that its optimal value is 0.001 among the candidate set of $\{1,0.1,0.01,0.001,0.0001\}$.} is used, and
\item $\text{o}_{\pi,k}:=o_{(2n+1)k+2n}$ is put through a softmax function to approximate $\pi_k(\mathbf{y})$. Specifically, we construct
\[ \hat{\pi}_k:= \frac{e^{\text{o}_{\pi,k}}}{\sum_{l=0}^{K-1} e^{\text{o}_{\pi,l}}}, \quad k=0,1,\cdots,K-1 \]
and their target values are $\pi_k(\mathbf{x})$. The point of using a softmax function is to ensure $\hat{\pi}_k\in[0,1]$ and $\sum_{k=0}^{K-1} \hat{\pi}_k = 1$.
\end{itemize}

Note that in practice the target values $\{\mu_k(\mathbf{y}),\sigma_k(\mathbf{y}),\pi_k(\mathbf{y})\}_{k=0}^{K-1}$ are not known. So when constructing the loss function MDN uses the method of Maximal-Likelihood. The idea of maximal likelihood estimation is used to set the loss for training the MDN. It is set as the negative of the estimated log likelihood of the data. Specifically,
\begin{gather}
\label{mdnLoss}
 \text{MDN Loss} =-\sum_{j=0}^{N-1} \log(\hat{p}_{\mathbfit{X}|\mathbfit{Y}=\mathbf{y}_j}(\mathbf{x}_j))= -\sum_{j=0}^{N-1} \log\left(\sum_{k=1}^K\hat{\pi}_k(\mathbf{y}_j)\hat{p}_k(\mathbf{x}_j)\right),
\end{gather}
where $\hat{p}_k$ denotes the pdf of $\mathcal{N}(\hat{\mu}_k(\mathbf{y}_j), \hat{\sigma}^2_k(\mathbf{y}_j))$. Therefore, minimizing the loss is equivalent to looking for weights in the FNN structure of MDN that maximize the log-likelihood of the dataset.

Based on above setting, at the end of this subsection, we construct a new model based on MDN. Our new model, named as FW-MDN, consists of an MDN with a modified loss function which consists of two parts. One is the ordinary MDN loss, and the other is the distance between $\hat{\mathbf{y}}$ and $\mathbf{y}$, where $\hat{\mathbf{y}}$ is the prediction made by the forward model using the MDN outputs. Specifically, our loss function is
\begin{gather}
\label{loss}
Loss := -\sum_{j=0}^{N-1} \log(\hat{p}_{\mathbfit{X|\mathbfit{Y}=\mathbf{y}_j}}(\mathbf{x}_j)) + \sum_{j=0}^{N-1} |\hat{\mathbf{y}}_j-\mathbf{y}_j|^2 + \alpha_b\sum_{l=0}^{B-1}b_l^2 + \alpha_w\sum_{l=0}^{W-1}w_l^2,
\end{gather}
where the first negative sum is the MDN loss defined in (\ref{mdnLoss}). In the second sum, $\hat{\mathbf{y}}_j = \sum_{k=0}^{K-1} \pi_k \hat{f}(\hat{\mu}_{k}(\mathbf{y}_j))$, where $\{\hat{\mu}_{k}\}$ are part of the MDN outputs (see, e.g., Figure \ref{MDN}) and $\hat{f}$ is the approximated\footnote{The forward model is approximated with an FNN using the training and validation sample sets. Its training $R^2$ is $96.5\%$ and validation $R^2$ is $96.4\%$.} forward model used for data generation. The third and fourth sum are regularizations for the bias terms $\{b_l\}_{l=0}^{B-1}$ and weights $\{w_l\}_{l=0}^{W-1}$ associated with the input, hidden, as well as the MDN output layers. Also, $\alpha_b$ and $\alpha_w$ are the regularization coefficients which will be treated as hyper-parameters.

\begin{figure}
	\centering
		\includegraphics[scale=.5]{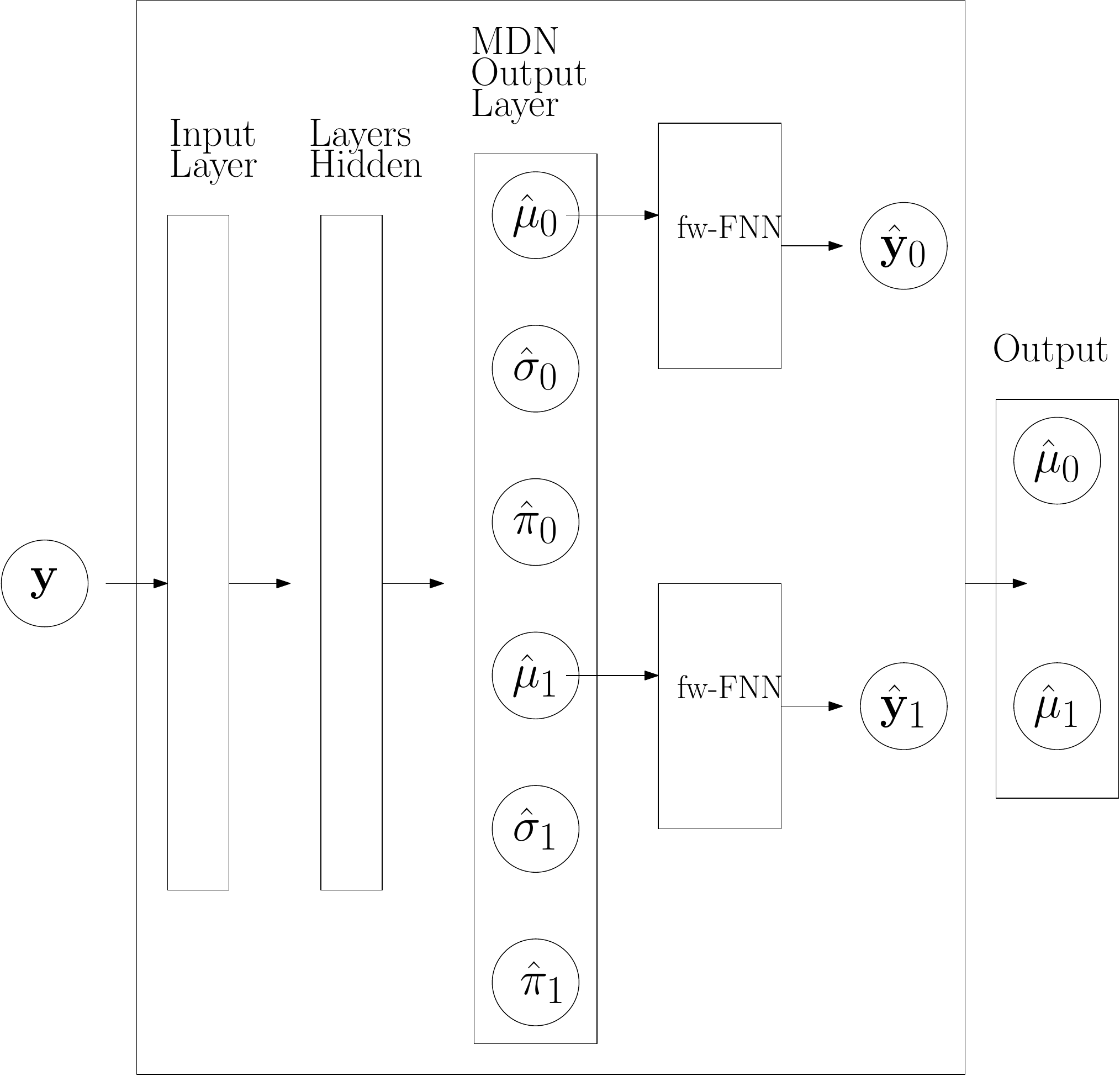}\\
       	\caption{Illustration of the FW-MDN Structure}
	\label{toyfig3}
\label{fwMDN}
\end{figure}

For each sample input $\mathbf{y}_j$, our FW-MDN outputs $\{\hat{\mu}_k\}_{k=0}^{K-1}$, the estimated means of the $K$ components of the MDN, each of which is supposed to track a branch of the target multi-valued function. \footnote{For each $\mathbf{y}_i$, the MDN outputs a set of parameters that determine a Gaussian mixture distribution, $\{\hat{\mu}_k(\mathbf{y}_i),\hat{\sigma}_k(\mathbf{y}_i),\hat{\pi}_k(\mathbf{y}_i) \}_{k=0}^{K-1}$. While our FW-MDN only outputs the K means.}
An example of FW-MDN for the case of $K=2$ and $\text{dim}(\mathbf{x})=1$ is illustrated in Figure \ref{fwMDN}, where fw-FNN denotes the pre-trained FNN that approximates the forward model, which is $\hat{f}$ in the loss. A toy example with MDN is provided in Appendix B.

\subsection{The training and testing of FW-MDN}
\label{fwMDN-pred}
Training FW-MDN is in fact training the involved FNN\footnote{Only the weights associated with the `Input Layer' and `Hidden Layers' in Figure \ref{fwMDN} are trained when training FW-MDN.}. Note that fw-FNN is pre-trained and its weights do not change when training FW-MDN. The model hyper-parameters are selected using the well known validation-set approach, which are shown in Table \ref{table1}.

\begin{table}[!htbp]
\caption{Model performances with different combinations of hyper-parameters and activations. $\alpha_b$ is the regularization coefficient for bias terms and $\alpha_w$ is the regularization coefficient for weights.}
\begin{center}
\begin{tabular}{p{1.5cm} p{4.2cm}  p{2.0cm} p{1.5cm} p{1.5cm} p{0.5cm} }
\toprule
Model & Hidden Layers and Nodes  &Activation  & $\alpha_b$ & $\alpha_w$ & $K$ \\
\toprule
fw-FNN     & $(40,100,200,200)$    &tanh                  &0.001   &0.001   & NA   \\
FW-MDN & $(400,300,300,300,300)$    &tanh      &0.00001   &0.00001  & 2   \\
\bottomrule
\end{tabular}
\label{table1}
\end{center}
\end{table}

Since each sample $(\mathbf{x}_j,\mathbf{y}_j)$ at most lies on one branch of the multi-valued function, it can only match one element in the FW-MDN output set $\{ \hat{\mu}_k\}_{k=0}^{K-1}$. Hence, we construct the measure for model performances as the distance between $\mathbf{x}_i$ and its closest neighbor in $\{ \hat{\mu}_k(\mathbf{y}_j)\}_{k=0}^{K-1}$. That is, we define the performance measure statistic as
\[ \mathcal{M} = 1 - \frac{\sum_{i=0}^{N-1}(\hat{\mu}^*_j - \mathbf{x}_j)^2}{\sum_{j=0}^{N-1}(\mathbf{x}_j - \bar{\mathbf{x}})^2}, \]
where $\bar{\mathbf{x}}:=\frac{1}{N}\sum_{j=0}^{N-1}\mathbf{x}_j$ and $\hat{u}^*_j:=\text{argmin}_{\hat{\mu}_k(\mathbf{y}_j),k\in \{0,1,2,...,K-1\}} | \hat{\mu}_k (\mathbf{y}_j) - \mathbf{x}_j |^2$. The construction of $\mathcal{M}$ is similar to that of the commonly used measure for model performance $R^2$, and the closer $\mathcal{M}$ is to 1 the better performance the FW-MDN has. This measure reflects FW-MDN's capacity to include the true $\mathbf{x}_j$ in its outputs, which is appropriate for any model that approximates a multi-valued function. We call $\hat{\mu}^*$ as the FW-MDN prediction for convenience, although the FW-MDN in fact predicts a set ${\{\hat{\mu}_k \}}_{k=0}^{K-1}$.

\begin{figure}[!htbp]
	\centering
		\includegraphics[scale=0.70]{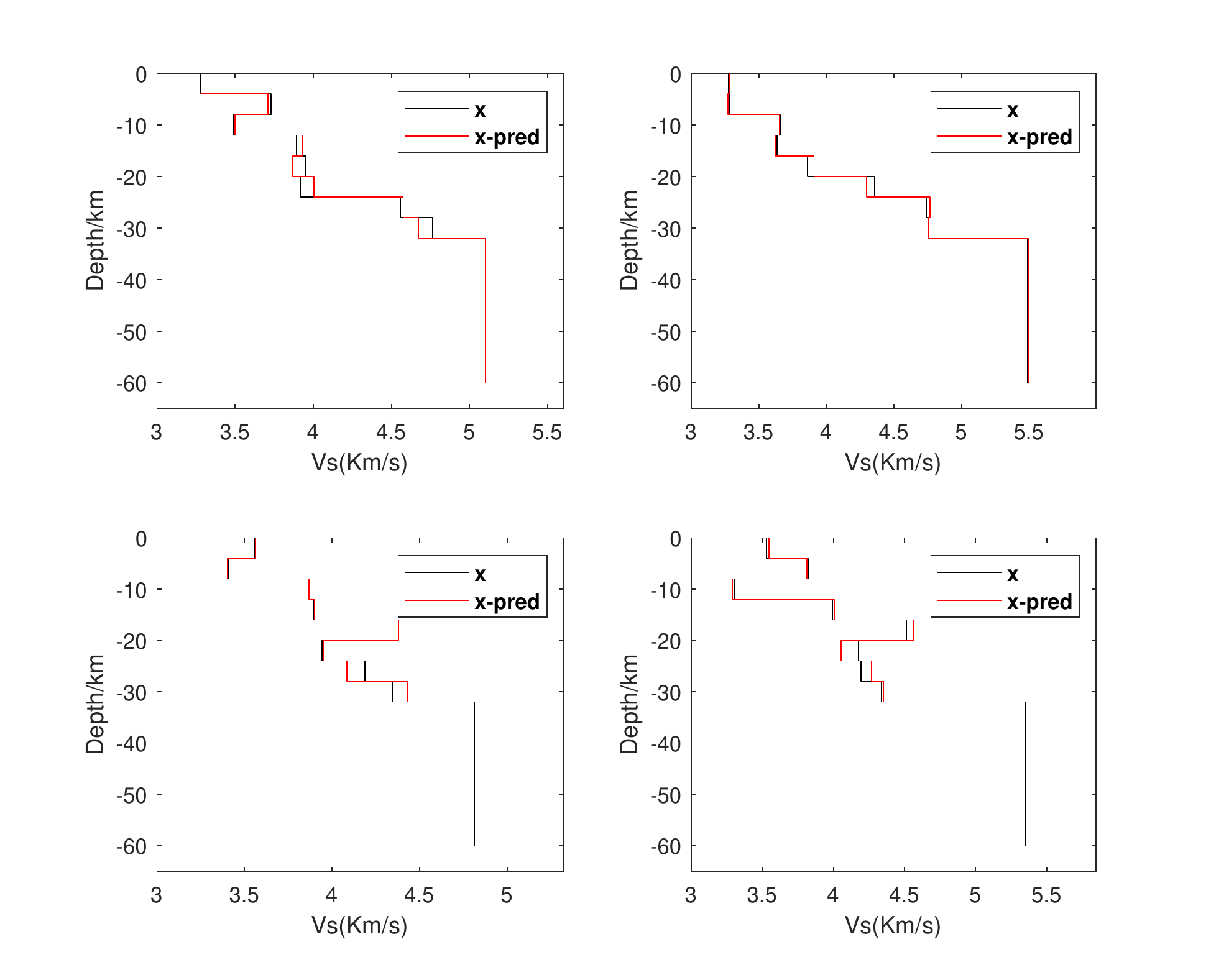}\\
       	\caption{The prediction $\hat{\mu}^*$ and and the true $\bar{V}_s$ (i.e., $\mathbf{x}$) for four randomly chosen samples}
\label{predictionPlot}
\end{figure}

\begin{figure}[!htbp]
	\centering
		\includegraphics[scale=0.5]{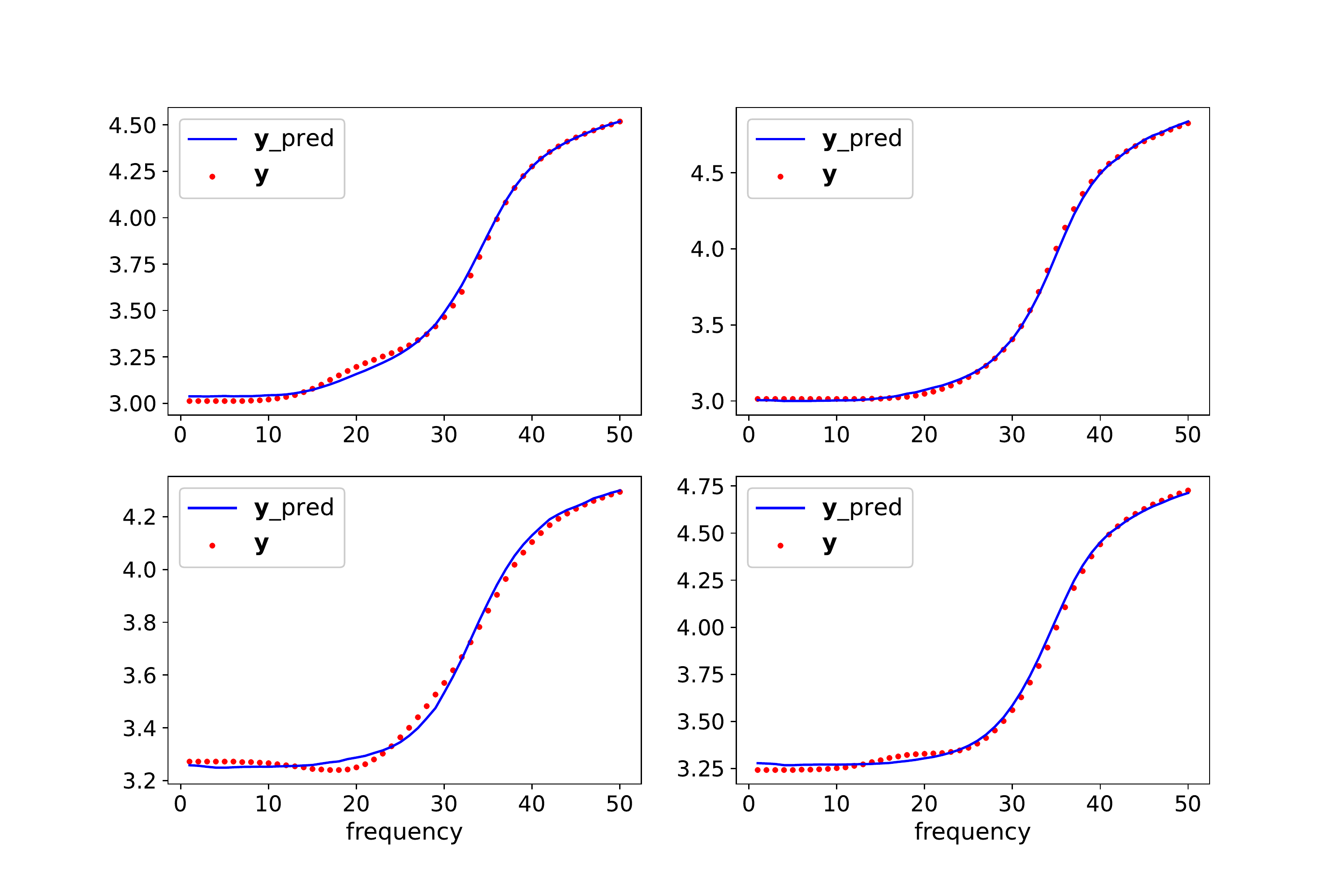}\\
       	\caption{Predicted dispersion curves ($\hat{\mathbf{y}}$) and the true values ($\mathbf{y}$)  for four randomly chosen samples}
\label{predictionyPlot}
\end{figure}

The testing performances of the trained FW-MDN are reported in Table \ref{testM}. Comparing to the ordinary FNN whose testing $R^2$ are reported in Table \ref{FNN9layer}, there is a $25\%$ increase in FW-MDN's overall performance, and its scores in the middle layers, say $\{x_5,x_6,x_7\}$, also improves significantly from below $40\%$ to above $80\%$. The predictions $\hat{\mu}^*$ for four randomly chosen samples are plotted\footnote{Figure \ref{predictionPlot} and Figure \ref{noisedPredPlot} are plotted with Matlab, while Figure \ref{predictionyPlot} and Figure \ref{noisedyPred} are plotted with Python.} in Fig. \ref{predictionPlot}.

\begin{table}[!htbp]
\caption{FW-MDN Test Performances.}
\begin{center}
\begin{tabular}{ p{1.5cm} p{1.5cm}  p{0.5cm} p{1.5cm} p{1.5cm}  }
Entry &   Test $R^2$  &        &   Entry        &  Test $R^2$         \\
\toprule
$x_0$	&	$99.7\%$	&		&	$x_5$	&	$87.5\%$		\\
$x_1$	&	$99.8\%$	&		&	$x_6$	&	$80.1\%$		\\
$x_2$	&	$99.1\%$	&		&	$x_7$	&	$85.3\%$		\\
$x_3$	&	$94.3\%$	&		&	$x_8$	&	$99.6\%$		\\
$x_4$	&	$84.1\%$	&		&	Overall&   $92.2\%$           \\
\bottomrule
\end{tabular}
\label{testM}
\end{center}
\end{table}

To see whether the FW-MDN prediction $\hat{\mu}^*$ produces dispersion curves close to the true values (i.e., the true $\hat{\mathbf{y}}$), on the test sample set we compute the $R^2$ between $\{\hat{f}(\hat{\mu}_j^*) \}_{j=0}^{N-1}$ and $\{\mathbf{y}_j\}_{j=0}^{N-1}$, and obtain a satisfying result of $94.8\%$, indicating that our FW-MDN prediction and the true value of $\bar{V}_s$ produce dispersion curves that are close to each other. Fig. \ref{predictionyPlot} plots $\hat{f}(\hat{\mu}^*)$ and $\mathbf{y}$ for the four samples in Fig. \ref{predictionPlot}.

\subsection{Robustness of FW-MDN under artificial noises}
As measurement errors in the dispersion curves are common, the model's robustness against noises needs to be investigated. Hence, we test the previously trained FW-MDN's performance on artificially noied testing samples.

Specifically, for each $(\mathbf{x}_j, \mathbf{y}_j)$ in the test set, define
\begin{gather}
\mathbf{y}^{\text{noised}}_j := \mathbf{y}_j(1+\varepsilon_j),
\label{uniform-err}
\end{gather}
where $\{\varepsilon_i\}$ are uniformly (and independently) drawn from $(-0.5\%, 0.5\%)$. Then we input $\{\mathbf{y}^{\text{noised}}_j\}$ to the FW-MDN and compute $\mathcal{M}$ between $\{ \mathbf{x}_j \}$ and their nearest neighbors $\{ \hat{\mu}_j^* \}$ in the FW-MDN outputs. The result is $ -19.6\% $, which indicates the previously trained FW-MDN is non-robust to measurement errors \footnote{If the noises drawn from $ N(0, 0.2\%^2) $ are added to the test data, FW-MDN's test performance becomes $18.7\%$, which also implies non-robustness to noises.}.

To create models that are robust to noises, we build new FW-MDNs \footnote{We keep model structures and hyper-parameters the same as those in Table \ref{table1}.} with the training and validation samples artificially noised with uniformly-distributed errors (the same way as the formula (\ref{uniform-err})), which enables the model to learn to discern noises. Specifically, when training the model we input the noised data to the model, but the loss (\ref{loss}) are computed with model outputs $\{\hat{\mathbf{y}}_j\}$ and un-noised data $\{(\mathbf{x}_j, \mathbf{y}_j)\}$.

\begin{table}[!b]
\caption{FW-MDN Performances $\mathcal{M}$ on Noised Test Samples.}
Each column records the performance of the FW-MDN on the test data noised with a error type. For example, in the column of $N(1\%$, $0.1\%^2)$, we independently draw $\varepsilon_i \sim N(1\%$, $0.1\%^2)$ and noise the test data in the way of (\ref{uniform-err}). Then we compute $\mathcal{M}$ between the test $\{\mathbf{x}_i\}$ and the FW-MDN outputs. The row of $\mathbf{x}$ reports the overall $\mathcal{M}$. The row $x_l, l=0,1,\cdots,8$ reports $\mathcal{M}$ for the $l^{th}$ entry of $\mathbf{x}$. The row of $\mathbf{y}$ reports the $R^2$ between $\mathbf{y}$ and $\hat{f}(\hat{\mu}^*)$. The row of $\mathbf{y}^{noised}$ reports the $R^2$ between $\mathbf{y}^{noised}$ and $\hat{f}(\hat{\mu}^*)$.
\begin{center}
\begin{tabular}{ p{0.9cm} p{1.5cm}  p{1.9cm} p{1.9cm} p{2.5cm} p{2.0cm} }
\toprule
 \multicolumn{6}{c}{Noise added to training and validation: Unif$(-0.8\%, 0.8\%)$} \\
\hline
Target       &Noise-free     &$N(0,$ $0.25\%^2)$  &$N(0,$ $0.50\%^2$)  &Unif$(-0.50\%,$ $0.50\%$)       &  Unif$(-0.90\%$, $0.90\%)$         \\
\toprule
$\mathbf{x}$	&	$76.0\%$	& $74.3\%$		&	$70.3\%$	&	$74.0\%$	& $70.0\%$	\\
$\mathbf{y}$& $94.7\%$  &$97.6\%$       &   $94.5\%$    &   $94.6\%$    & $94.5\%$\\
$\mathbf{y}^{noised}$& NA  &$94.3\%$       &   $93.3\%$    &   $94.2\%$    & $93.2\%$\\
$x_0$	&	$99.7\%$	& $99.7\%$		&	$99.6\%$	&	$99.7\%$    & $99.6\%$ 		\\
$x_1$	&	$99.3\%$	& $98.9\%$		&	$97.8\%$	&	$98.7\%$	& $97.7\%$	\\
$x_2$	&	$92.2\%$	& $89.7\%$		&	$83.4\%$	&	$89.1\%$	& $82.8\%$	\\
$x_3$	&	$73.0\%$	& $69.5\%$		&	$62.5\%$	&	$69.2\%$	& $61.6\%$	\\
$x_4$	&	$71.8\%$	& $69.4\%$		&	$63.6\%$	&   $68.7\%$    & $64.0\%$   \\
$x_5$	&	$58.9\%$	& $58.0\%$		&	$55.3\%$	&	$57.9\%$	& $54.6\%$	\\
$x_6$	&	$49.7\%$	& $47.8\%$		&	$42.3\%$	&	$47.8\%$	& $43.1\%$	\\
$x_7$	&	$39.1\%$	& $36.3\%$		&	$28.5\%$	&	$35.0\%$	& $26.6\%$	\\
$x_8$	&	$99.9\%$	& $99.9\%$		&	$99.7\%$	&	$99.9\%$	& $99.7\%$	\\

\bottomrule
\end{tabular}
\label{noised-test}
\end{center}
\end{table}

\begin{figure}[!htb]
	\centering
		\includegraphics[scale=0.70]{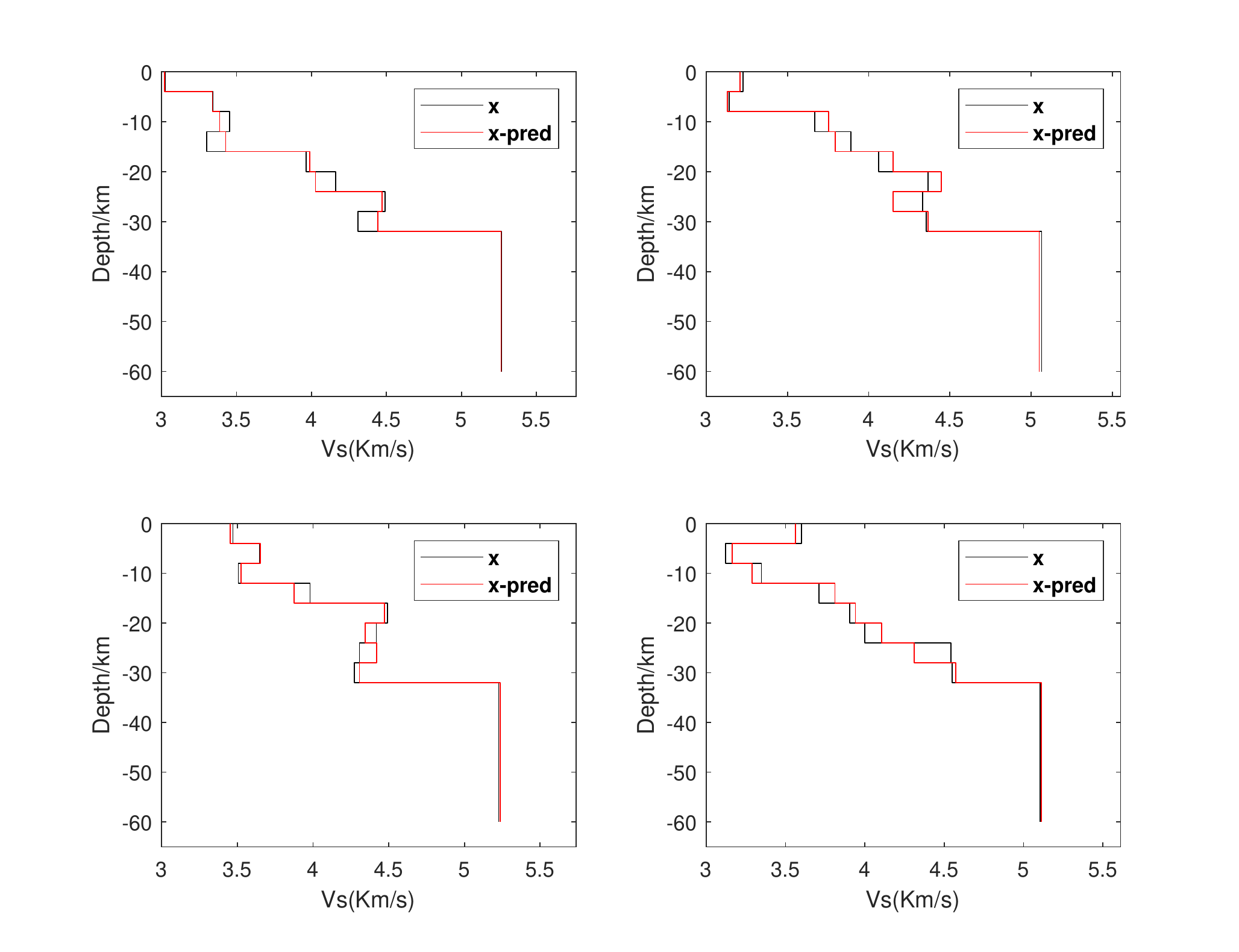}\\
       	\caption{The prediction $\hat{\mu}^*$ and the true $\bar{V}_s$ (i.e., $\mathbf{x}$) for four randomly chosen noised test samples}
\label{noisedPredPlot}
\end{figure}

\begin{figure}[!htb]
	\centering
		\includegraphics[scale=0.45]{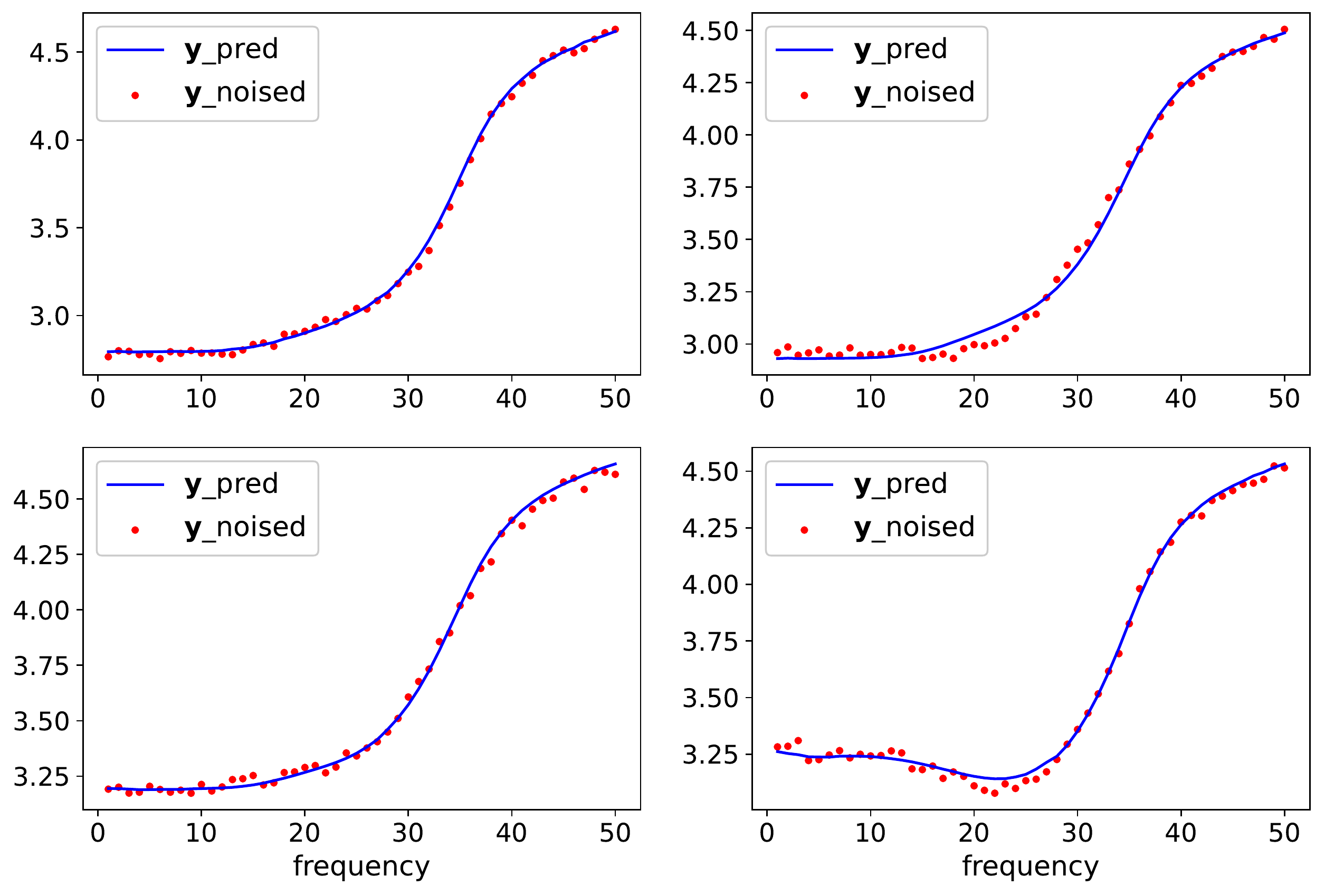}\\
       	\caption{Plots of $ \mathbf{y}^{noised}$ and $\hat{f}(\hat{\mu}^*)$  for the four samples in Figure \ref{noisedPredPlot}}
\label{noisedyPred}
\end{figure}

To determine how much noises should be added to the training dataset, for each of the noise type $\{\text{Unif}(-0.2\%, 0.2\%)$, $\text{Unif}(-0.5\%, 0.5\%)$, $\text{Unif}(-0.8\%, 0.8\%)$, \\
$\text{Unif}(-1\%, 1\%)\}$ we train a FW-MDN, and select the one with the best performance on the validation data noised with errors uniformly distributed in $(-1\%,1\%)$. The selected noise is $\text{Unif}(-0.8\%, 0.8\%)$, with which we noise the training and validation data and build a new FW-MDN. In Table \ref{noised-test} we report its performance on testing samples noised with various types of artificial errors. The results show that, comparing with the FW-MDN trained with noise-free data, the newly-trained FW-MDN is much more robust to noises. For example, the overall $\mathcal{M}$ on the test samples with various type of noises are all above $70\%$, which are satisfactory in most real-world cases. Figure \ref{noisedPredPlot} plots $\hat{\mu}^*$ obtained from the model and $\mathbf{x}$ for four randomly drawn testing samples that are noised with Unif$(-0.9\%,0.9\%)$ errors, from which we can see that the predictions $\hat{\mu}^*$ are close to the true values of $\mathbf{x}$. Figure \ref{noisedyPred} plots $\mathbf{y}^{noised}$ and $\hat{f}(\hat{\mu}^*)$ for the four samples in Figure \ref{noisedPredPlot}, in which the noised dispersion curves $\mathbf{y}$ are also close to the predictions $\hat{f}(\hat{\mu}^*)$ obtained by forwardly propagating the FW-MDN outputs. All these results show that the newly-trained FW-MDN is robust to noises.


\section{Conclusion}
\label{Conclusion}
In this work, we develop a new Rayleigh wave dispersion equation (\ref{f}) (or (\ref{fjVector})), from which the dispersion curve can be easily computed. Based on this new forward model, we use the structure of a mixture density network to develop a new inversion solver, named FW-MDN, whose training loss function involves the forward model, to reconstruct the S-wave velocity from the data of dispersion curves. The FW-MDN is capable to deal with the non-uniqueness phenomenon appeared in our inverse geophysics problem, which other existing studies that apply machine learning fail to accommodate.

Although more experiences with real-word data is needed to fully understand the potential and the limitations of FW-MDN, the initial applications to the model problem and the artificial data system are promising, yielding excellent fits of the data and well-defined S-wave velocity. Moreover, the proposed novel forward model and the machine learning based inversion solver can be easily extended to the case with Rayleigh wave multiple mode dispersion curve. We therefore believe FW-MDN will be a useful tool to study many inverse problems in geophysics.

\bibliography{refpool}  
\bibliographystyle{plain}

\section*{Appendix A: performances of ordinary Feedforward Neural Networks (FNNs) and its failure in multi-layer  medium model}

In this appendix, we approach the inversion problems (\ref{fjVector}) with ordinary FNNs. For each case we build an ordinary FNN to predict the S-wave velocity $\bar{V}_S$ from the dispersion curve $\mathbf{C}$. The hyper-parameters, such as the number of hidden layers and nodes, regularization coefficients, and activation functions, are selected using the validation-set approach. The performance of the optimal FNN in each case are listed in Table \ref{FNN-performances}, which shows that FNN can solve the inversion problem perfectly for the 3-layer and 5-layer case. However, its performance decreases sharply for the 9-layer case. Table \ref{FNN9layer} reports the FNN's performance on each entry of the input vector, we can see that the predictions for the middle layer $\bar{V}_S$ are off. For example, the $R^2$ is merely $30.7\%$ for Vs of in the 8th layer ($x_7$).

\begin{table}[!htbp]
\caption{Model performances with different combinations of hyper-parameters and activations. $R^2 = 1 - \frac{\sum_{j=0}^{N-1}(\hat{\mathbf{x}}_j - \mathbf{x}_j)^2}{\sum_{i=0}^{N-1}(\mathbf{x}_j - \bar{\mathbf{x}})^2}$, where $\bar{\mathbf{x}}:=\frac{1}{N}\sum_{j=0}^{N-1}\mathbf{x}_j$. $\alpha_b$ is the regularization coefficient for bias terms and $\alpha_w$ is the regularization coefficient for weights.}
\begin{center}
\begin{tabular}{p{1.2cm} p{4.5cm}  p{1.5cm} p{1.0cm} p{1.0cm} p{0.9cm} p{0.9cm}}
\toprule
case&Hidden Layers and Nodes  &Activation  & $\alpha_b$ & $\alpha_w$ & Train $R^2$ &Test $R^2$ \\
\toprule
3-layer &$(150,150,150,100,100)$    &tanh      &0.001   &0.001  &$99.7\%$  &$99.7\%$   \\
5-layer &$(150,150,150,100,100)$    &tanh      &0.001   &0.001  &$98.7\%$  &$98.7\%$   \\
9-layer &$(400,400,300,300,300)$    &tanh      &0.001   &0.001  &$66.5\%$  &$66.4\%$   \\
\bottomrule
\end{tabular}
\label{FNN-performances}
\end{center}
\end{table}

\begin{table}[!htbp]
\caption{FW-MDN Test Performances.}
\begin{center}
\begin{tabular}{ p{2cm} p{2cm}  p{2cm} p{2cm} p{2cm}  }
\toprule
Entry &   Test $R^2$  &        &   Entry        &  Test $R^2$         \\
\toprule
$x_0$	&	$98.3\%$	&		&	$x_5$	&	$29.5\%$		\\
$x_1$	&	$96.7\%$	&		&	$x_6$	&	$37.2\%$		\\
$x_2$	&	$81.2\%$	&		&	$x_7$	&	$30.7\%$		\\
$x_3$	&	$48.2\%$	&		&	$x_8$	&	$99.3\%$		\\
$x_4$	&	$37.4\%$	&		&	Overall	&	$66.5\%$           \\
\bottomrule
\end{tabular}
\label{FNN9layer}
\end{center}
\end{table}

\begin{figure}[!htbp]
	\centering
		\includegraphics[scale=.45]{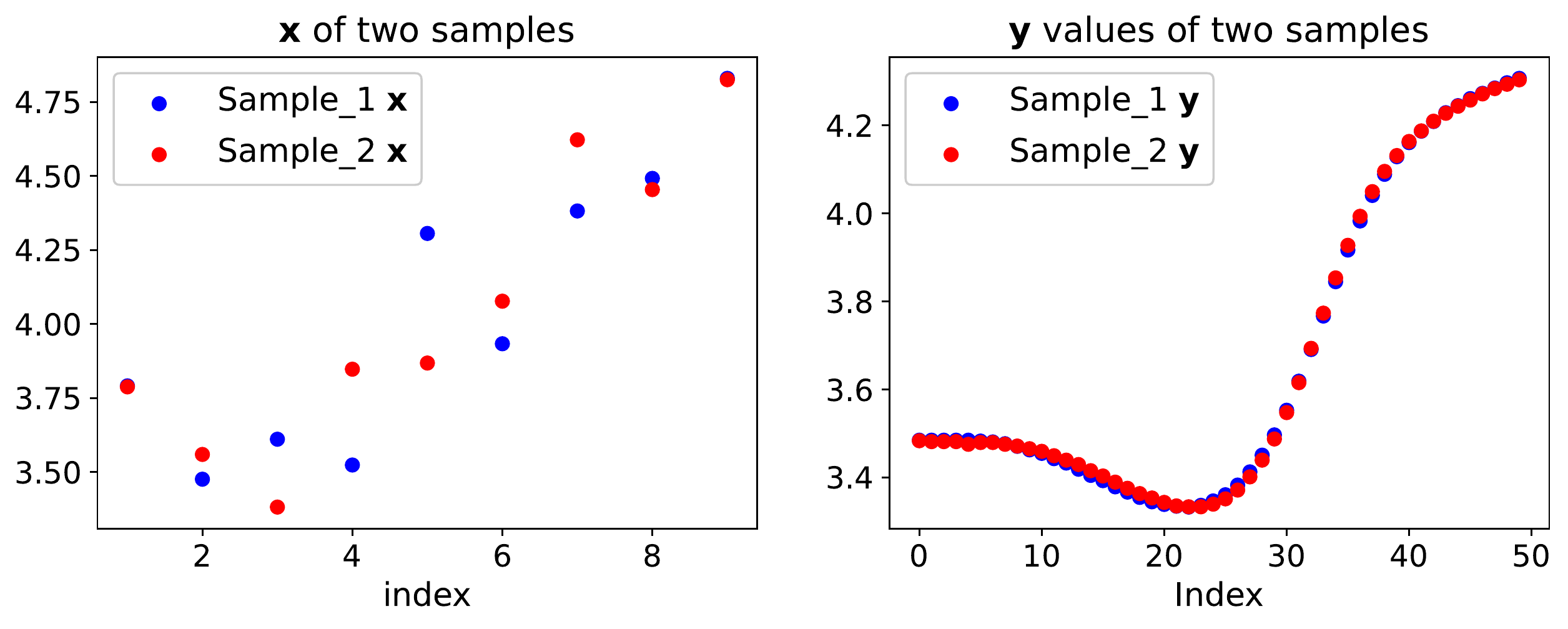}
       	\caption{Plots of 200 samples in the artificial dataset of our problem.}
	\label{onetomany}
\end{figure}

Now, let us discuss the reasons for the FNN under-performance for the case of 9-th layer medium model of the Earth. Actually, among our simulated 9-layer data, we find that certain data samples have very close $\mathbf{y}$ values, yet their $\mathbf{x}$ values differ significantly. Figure \ref{onetomany} illustrates this phenomenon with two samples from our dataset. This is the well-known non-uniqueness problem. It implies that the target function for our learning model is a multi-valued function\footnote{Suppose $f$ is an ``non-injective" mapping from a set $\mathcal{X}$ to another set $\mathcal{Y}$, i.e., for certain $\mathbf{y}\in \mathcal{Y}$ there are $\{\mathbf{x}_k\}_{k=0}^{K-1} \subset \mathcal{Y}$ such that $f(\mathbf{x}_k)=\mathbf{y}$. Then the inverse mapping, denoted by $f^{-1}$, is called a \textit{multivalued} function. If there also exist $K$ single-valued functions $\{f_k^{-1}(\mathbf{y})\}_{k=0}^{K-1}$ such that $f^{-1}_k(\mathbf{y})=\mathbf{x}_k$, then each $g_k$ is called a \textit{branch} of $f^{-1}$. A multi-valued function in fact refers to a relation $F(\mathbf{x},\mathbf{y})=0$ such that for a certain value of $\mathbf{y}$, say $\mathbf{y}_0$, there are more than one values of  $\mathbf{x}$, say $\{\mathbf{x}_0,\mathbf{x}_1\}$, such that $F(\mathbf{x}_k,\mathbf{y}_0)=0$ for $k=0,1$.} which an ordinary FNN is not able to approximate. This is because an FNN is a single-valued function. A recognized learning model that can handle this approximation problem is the mixture density network, which is introduced in Sect \ref{introMDN}.

\section*{Appendix B: a toy example with MDN}

In this appendix, we illustrate the performance of an MDN with a toy example, where the dataset is generated based on the equation $y = x+0.3\sin(2\pi x)$ where $x\in[0,1]$. Specifically,
\begin{itemize}
\item We generate $10,000$ values from the uniform distribution Uniform(0,1) and denote them by $\{x_i\}_{i=0}^{9999}$.
\item We generate $10,000$ values from the uniform distribution Uniform$(-0.1,0.1)$ and denote them by $\{\epsilon_i\}_{i=0}^{9999}$.
\item For $i\in \{0,\cdots,9999\}$, define $y_i = x_i  + 0.3\sin(2\pi x_i) +  \epsilon_i$, where $\{\epsilon_i\}$ serves as noises.
\end{itemize}

\begin{figure}[h]
	\centering
		\includegraphics[scale=0.55]{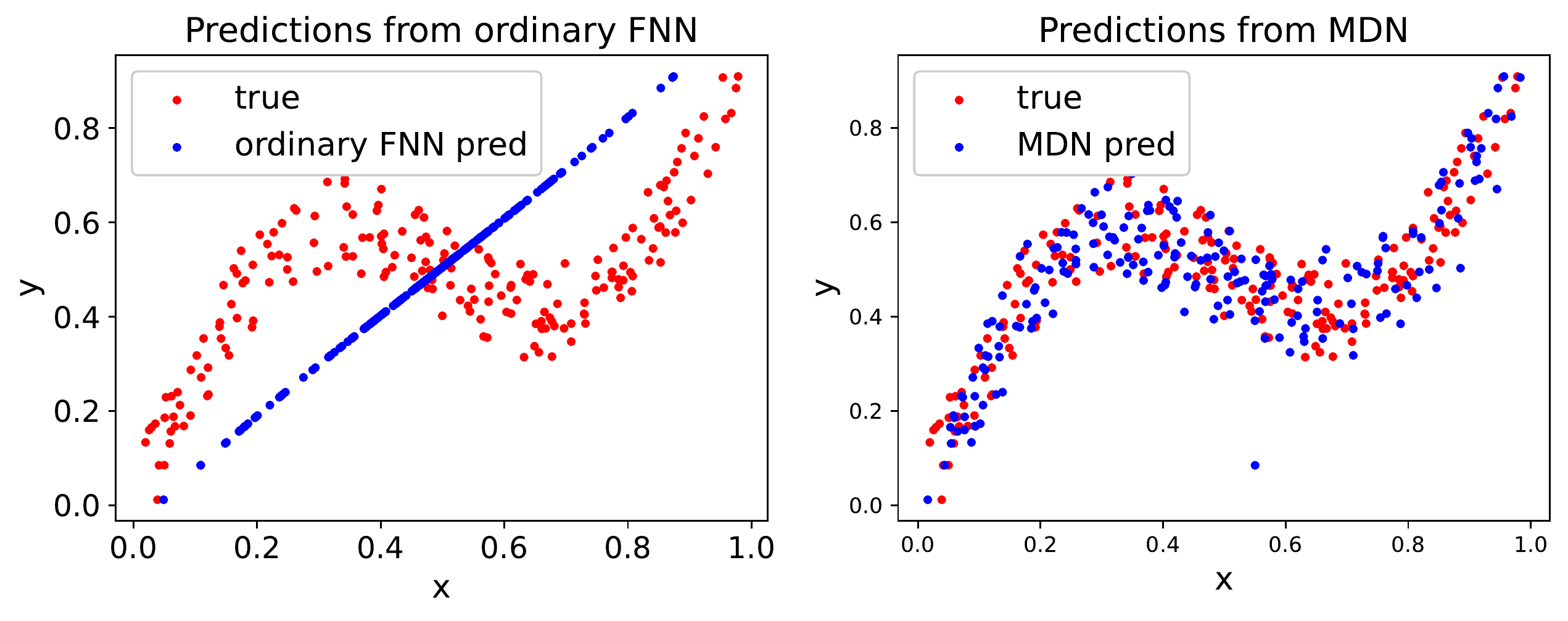}\\
       	\caption{Toy Example: Plots of 200 FNN predictions}
	\label{toyfig}
\end{figure}

The plots of 200 samples, $\{(x_i, y_i)\}_{i=0}^{199}$, are the red dots in Figure \ref{toyfig}, which shows that one value of $y$ may correspond to more than one values of $x$.

To show the advantage of MDN, we also report the performance of an ordinary FNN on this toy example, which intakes $y_i$ to predict $x_i$. The training \footnote{We split the 10,000 data to training, validation, and testing sets, and use the validation-set approach to select hyper-parameters like number of hidden layers, activation functions, and regularization coefficients.} $R^2$ of the ordinary FNN is merely $37.6\%$, which is expected since FNN is a continuous function, and for one value of $y$ it outputs only one value of $x$, which makes it incapable of finding more than one values of $x$ corresponding to any particular value of $y$. The predictions corresponding to the 200 samples are plotted on the left of Figure \ref{toyfig}. Clearly, the predictions are inconsistent to the true values.

We then switch to a MDN model with a mixture of four ($K=4$) Gaussian distributions. After training the model, we randomly pick 200 samples $\{x_i\}_{i=0}^{199}$ from the test data. For each $y_i$ we input it to the trained MDN\footnote{In fact, we input $y_i$ to the trained FNN involved in the MDN and get the outputs $\{\hat{\pi}_k(y_i),\hat{\mu}_k(y_i), \hat{\sigma}_k(y_i)\}_{k=1}^4$.}. Then we draw a sample, denoted by $\hat{x}_i$, from the Gaussian mixture with parameters $\left\{\pi_k(y_i),\mu_k(y_i), \sigma_k(y_i) \right\}_{k=1}^4$. The plots are on the right of \ref{toyfig}. You can see that the spreading of these 200 predictions are close to the 200 true values. 

\end{document}